\documentclass[a4paper]{article}

\usepackage{lineno}
\usepackage{amsmath,amsfonts,amssymb}
\usepackage{graphicx}
\usepackage{tikz}
\usepackage{float}
\newcommand{\D}{\displaystyle}

\newcommand{\veps}{\varepsilon}
\newcommand{\dsum}{\D\sum}
\newcommand{\dint}{\D\int}

\newcommand{\muf}{\mu_{ph}}	
\newcommand{\mug}{\mu_g}

\newcommand{\uk}{\textbf{u}}

\newcommand{\vk}{\textbf{v}}
\newcommand{\wk}{\textbf{w}}
\newcommand{\zk}{\textbf{z}}

\newcommand{\unmu}{U_N(\muk)}

\newcommand{\fk}{\textbf{f}}

\newcommand{\intOo}[1]{\dint_{\Omega_o(\mug)} #1\; d\Omega_o}
\newcommand{\intOr}[1]{\dint_{\Omega_r} #1\; d\Omega_r}
\newcommand{\intK}[1]{\dsum_{K\in \mathcal{T}_h}\dint_K #1\;d\Omega}

\newcommand{\dOr}{d\Omega_r}
\newcommand{\Or}{\Omega_r}
\newcommand{\cD}{\mathcal{D}}
\newcommand{\R}{\mathbb{R}}
\newcommand{\Tnmu}{T_N^\mu}

\newcommand{\TN}[1]{T_N(#1;\mu)}

\newcommand{\ra}{\rightarrow}

\newcommand{\normlur}[1]{\|#1\|_{0,1,\R}}
\newcommand{\normlt}[1]{\|#1\|_{0,3,\Omega}}
\newcommand{\normlc}[1]{\|#1\|_{0,4,\Omega}}
\newcommand{\normld}[1]{\|#1\|_{0,2,\Omega}}
\newcommand{\normlinf}[1]{\|#1\|_{0,\infty,\Omega}}
\newcommand{\normh}[1]{\|#1\|_{1,2,\Omega}}

\newcommand{\normX}[1]{\|#1\|_X}

\newcommand{\difZ}{Z^1_h-Z^2_h}

\newcommand{\cR}[1]{\mathcal{R}(#1;\muk)}
\newcommand{\DA}[1]{\mathcal{DA}(#1;\muk)}
\newcommand{\DAn}{\DA{\unmu}}
\newcommand{\Hw}[1]{H(#1;\muk)}
\newcommand{\difH}{\Hw{Z_h^1}-\Hw{Z_h^2}}
\newcommand{\en}{\epsilon_N(\muk)}
\newcommand{\enk}{\epsilon_N(\muk)}

\newcommand{\taunk}{\tau_N(\muk)}

\newcommand{\dnk}{\Delta_N(\muk)}
\newcommand{\bmu}{\beta_N(\muk)}
\newcommand{\dual}[1]{\left<#1,V_h\right>}
\newcommand{\ximu}{\xi(\muk)}

\newcommand{\xip}{\xi^p}
\newcommand{\zetav}{\zeta^\textbf{v}}
\newcommand{\phit}{\varphi^{\theta}}

\newcommand{\Th}{\mathcal{T}_h}
\newcommand{\muk}{\boldsymbol{\mu}}

\newcommand{\pih}{\Pi_h^*}

\newcommand{\dx}{\partial_x}
\newcommand{\dy}{\partial_y}
\newcommand{\red}[1]{\textcolor{black}{#1}}

\newcommand{\black}[1]{\textcolor{black}{#1}}

\newtheorem{lemma}{Lemma}
\newtheorem{proposition}{Proposition}
\newtheorem{theorem}{Theorem}
\newtheorem{Remark}{Remark}

\begin{document}

%\begin{frontmatter}
\title{Certified Reduced Basis VMS-Smagorinsky model for natural convection flow in a cavity with variable height}
\date{}

\author{
Francesco Ballarin\thanks{mathLab, Mathematics Area, SISSA, International School for Advanced Studies, via Bonomea 265, I-34136 Trieste, Italy. francesco.ballarin@sissa.it, gianluigi.rozza@sissa.it} \and 
Tom\'{a}s Chac\'{o}n Rebollo\thanks{IMUS \& Departamento de Ecuaciones Diferenciales y An\'{a}lisis Num\'{e}rico 
, Apdo. de correos 1160, Universidad de Sevilla, 41080 Seville, Spain. chacon@us.es, edelgado1@us.es}  
\and Enrique Delgado \'{A}vila\footnotemark[2]
\and Macarena G\'{o}mez M\'{a}rmol\thanks{Departamento de Ecuaciones Diferenciales y An\'{a}lisis Num\'{e}rico, Apdo. de correos 1160, Universidad de Sevilla, 41080 Seville, Spain. macarena@us.es}
\and 
Gianluigi Rozza\footnotemark[1]
}

\maketitle

\begin{abstract}
In this work we present a Reduced Basis VMS-Smagorinsky Boussinesq model, applied to natural convection problems in a variable height cavity, in which the buoyancy forces are involved. We take into account in this problem both physical and geometrical parametrizations, considering the Rayleigh number as a parameter, so as the height of the cavity. We perform an Empirical Interpolation Method to approximate the sub-grid eddy viscosity term that let us obtain an \textit{affine} decomposition with respect to the parameters. We construct an \textit{a posteriori} error estimator, based upon the Brezzi-Rappaz-Raviart theory, used in the greedy algorithm for the selection of the basis functions. Finally we present several numerical tests for different \red{parameter configuration}.
\end{abstract}


%\end{frontmatter}

%\linenumbers

\section{Introduction}

Nowadays, several industrial processes need numerical simulations, which are usually performed with the widespread high-fidelity approximation techniques such as finite element (FE), finite volumes or spectral methods, and they usually take very long time for computation \red{(several hours, even days)}. In many situations, the model that represents the behavior of an industrial process is given by a Partial Differential Equation (PDE)  depending on parameters. Reduced-order modeling (ROM) is used in parametrized PDE in order to try to reduce the high computational time required by its numerical solution, when large number of simulations with different parameter values are needed \cite{hestaven,holmes,Nguyen2005,Prudhomme2004,LibroBR,Stokes1,Rozza2008,Drohmann2015}.

In the context of fluids dynamics, even with the reduction of the computational cost provided by turbulent models, such as the Variational Multi-Scale (VMS) models
(\textit{cf.} \cite{LES-VMS}), with respect to the DNS, it is still expensive to compute accurately the real flows that commonly appear in industry problems, specially in cases where parameters play important roles. When a high number of computations for a fluid flow depending on parameters is required, ROM becomes useful. Several works for Reduced Basis (RB) model have been presented for Stokes equations \cite{Stokes2,Stokes1,Stokes3} and Navier-Stokes equations \cite{Deparis,Dparis-Rozza,Manzoni,Patera,Yano2014}. Most recently, there have been developed works for RB turbulent models, such as the Smagorinsky model \cite{PaperSmago}, or the VMS-Smagorinsky model \cite{Chacon2018}. In those last works, the ROM is constructed from the turbulent model. A different way to build ROM model for turbulent flows is the one in which initially we would build a ROM for Navier-Stokes equations, and then model the unresolved scales (either by eddy diffusion or other techniques) to build the turbulence model. This approach has been followed in \cite{iliescu,VMS-Iliescu}, for instance. We address here the construction of a RB model for the Boussinesq equations, that includes turbulent diffusivity for both momentum and energy equations. The turbulent diffusivities are modeled by a VMS-Smagorinsky approach, in such a way that eddy diffusion effects only act on the small resolved scales. One of the main advantages of using the VMS-Smagorinsky is that the eddy viscosity and eddy diffusivity only
affect the small resolved scales, avoiding over-diffusive effects, that can be maintained in the ROM setting thanks to the Reduced Basis framework that we have developed. \red{On the counterpart, we deal with non-linear terms related with the VMS-Smagorinsky model, that force us to use techniques such as the Empirical Interpolation Method, resulting a feasible ROM model.}

In this work we consider the application of a RB Boussinesq VMS-Smago-\\rinsky model to simulate a natural convection in a variable height cavity. In applications to architecture, this cavity represents a courtyard inside a building. The study of the heat exchange between the air and the walls inside the courtyard is of high interest to minimize the energy needs of the building.
The variability of the cavity height is considered through a geometrical parametrization of the domain. Since we are interested in solving efficiently the parameter-dependent problem, we need to reformulate the Boussinesq VMS-Smagorinsky model in a parameter-independent domain with a change of variables. With this change of variables, we obtain operators that depend on both physical and geometrical parameters. This setting, besides the Empirical Interpolation Method (EIM) (\textit{cf.} \cite{EIM2,EIM1}) for the non-linear eddy diffusivities, lets us approximate these non-linear operators by operators that depend affinely with respect to the parameters, both of physical and geometrical type. Then it is possible to store parameter-independent matrices and tensors in the offline phase. 

We tackle in this work some of the intermediate difficulties which should necessarily be solved when dealing with RB modelling of natural convection turbulent flows. Particularly, we analyze how to deal with the temperature when constructing the \textit{a posteriori} error estimator. From the numerical analysis point of view, we present the development of an \textit{a posteriori} error bound based upon the Brezzi-Rappaz-Raviart (BRR) theory \cite{BRR}, used in the snapshot selection in the greedy algorithm. This \textit{a posteriori} error estimator is an extension of the ones presented for the Navier-Stokes equations \cite{Deparis,Dparis-Rozza,Manzoni} and the Smagorinsky model \cite{PaperSmago}. The main difference for the \textit{a posteriori} error bound presented in this paper with the previous ones, is the necessity of considering a mollifier for the thermal eddy diffusivity term, due to the fact this term is no longer Lipschitz-continuous. Thanks to the consideration of this regularized term, we are able to develop the \textit{a posteriori} error estimator.

We present four different numerical tests for the buoyancy-driven cavity problem. In the first two tests, we consider a fixed height, i.e. we consider the geometrical parameter $\mu_g=1$, for different ranges of the Rayleigh number, ranging from moderate Rayleigh numbers values $Ra\in[10^4,10^5]$, to high Rayleigh numbers values $Ra\in[10^5,10^6]$. In the third one, we fix the Rayleigh number, with a moderate value $Ra=10^5$, and only the geometrical parameter changes. This test intends to represent a situation in which the environmental conditions are fixed and we are only interested in simulating the flow in cavities. In the last test, both the Rayleigh and the geometric parameter are taken into account. This test is more complex since two parameters are considered, thus the number of basis functions to include in our RB spaces increases with respect to the previous one. At the same time, the speed-up ratio decreases somewhat, although it remains on values around 50.

The paper is structured as follows: in section \ref{chap:Geom::sec:FE} we define the high fidelity problem in the reference domain from the one defined in the original domain depending on the geometric parameter. Then, in section \ref{chap:Geom::sec:RB}, we present the Reduced Basis problem, with the EIM approximation for the eddy viscosity and eddy diffusivity terms. In section \ref{chap:Geom::sec:a_post} we construct the \textit{a posteriori} error estimator. Finally in section \ref{chap:Geom::sec:Num}, we present the numerical results for the tests previously described, programmed in FreeFem++ (\textit{cf.} \cite{freefem++}). Conclusions are then summarized in section \ref{sec:Conclusions}.

\section{Continuous problem and full order discretization}\label{chap:Geom::sec:FE}

The aim of this work is to present a reduced order model for natural convection problems over domains with variable geometry. For this purpose, we propose a turbulence Smagorinsky model in which the buoyancy forces are modeled by the Boussinesq approach. \red{Let $\muk=(\muf,\mug)\in\R^{p_{ph}\times p_g}$ and} let $\Omega_o(\mu_g)$ \red{be} a bounded polyhedral domain in $\R^d$  $(d=2,3)$, depending on \red{the geometrical parameters}, commonly called original domain in the RB framework. Let $\Gamma(\mu_g)=\Gamma_D(\mu_g)\cup\Gamma_N(\mu_g)$ the Lipschitz-continuous boundary of $\Omega_o(\mu_g)$, where $\Gamma_D$ is the part of the boundary with Dirichlet conditions and $\Gamma_N$ the part of the boundary with Neumann conditions. 

We next present the continuous Boussinesq-Smagorinsky model that we consider in this work. Although the Smagorinsky approach is intrinsically discrete, we present it in a continuous form in order to clarify its relationship with the standard Boussinesq model:
\begin{equation}\label{chap:Geom::pb:Cont_Bouss}\left\{\begin{array}{ll}
\uk_o\cdot\nabla_o\uk_o-Pr\Delta_o\uk_o-\nabla_o\cdot(\nu_T(\uk_o)\nabla_o\uk_o)\\
+\nabla_o p_o-Pr\,\mu_{ph}\,\theta_o\,\textbf{e}_d=\fk&\mbox{ in }\Omega_o(\mug)\vspace{0.2cm}\\

\nabla_o\cdot\uk_o=0&\mbox{ in } \Omega_o(\mug)\vspace{0.2cm}\\

\uk_o\cdot \nabla_o \theta_o-\Delta_o\theta_o-\nabla_o\cdot(K_T(\uk_o)\nabla_o\theta_o)=Q&\mbox{ in } \Omega_o(\mug)\vspace{0.2cm}\\
\uk_o=0&\mbox{ on }\Gamma\vspace{0.2cm}\\

\theta_o=\theta_D&\mbox{ on }\Gamma_{D}\vspace{0.15cm}\\
\partial_n\theta_o=0&\mbox{ on }\Gamma_N\vspace{0.15cm}\\
\end{array}\right.
\end{equation} 

Here, $\uk_o$ is the velocity field, $p_o$ is the pressure and $\theta_o$ is the temperature. In addition, $\textbf{e}_d$ is the last vector of the canonical basis of $\R^d$, while $\muf$ and $Pr$ are the Rayleigh and Prandtl dimensionless numbers respectively. %Let us denote $\muk=(\muf,\mug)\in\cD$ the parameters considered. 
Both the external body forces $\textbf{f}$, and the heat source term $Q$, are given data for the problem. In (\ref{chap:Geom::pb:Cont_Bouss}), $\nu_T(\uk)$ is the eddy viscosity term, and $K_T(\uk)$ is the eddy diffusivity given by
\begin{equation}\label{chap:Bouss::eq:eddy_conduct}
K_T(\uk)=\dfrac{1}{Pr}\nu_T(\uk).
\end{equation}

In (\ref{chap:Geom::pb:Cont_Bouss}), we represent by $\theta_D$ a given temperature over the boundary $\Gamma_D$. For simplicity of the analysis, we further consider that $\theta_D=0$. In the case of considering non-homogeneous boundary conditions, it is enough to define a lift function $\theta_g$ such that $\theta_g|_{\Gamma_D}=\theta_D$. In \cite{PaperSmago}, an analysis with the lift function is already done for a RBM Smagorinsky model.

Let us consider the spaces $Y^o=(H_0^1(\Omega_o))^d$, $M^o=L_0^2(\Omega_o)$, $\Theta^o=H_0^1(\Omega_o)$. We denote by $\|\cdot\|_{k,p,\Omega}$ the norm of the Sobolev space $W^{k,p}(\Omega)$. We consider the $H^1$-seminorm for the velocity and temperature spaces, and the $L^2$-norm for the pressure space. In addition, let us define the Sobolev embedding constants $C_u$ and $C_\theta$, associated to these norms, such that
\begin{equation}\label{chap:Geom::eq:Sobolev_u}
\normlc{\vk}\le C_u \normld{\nabla \vk}, \quad\forall \vk\in Y,
\end{equation}
and
\begin{equation}\label{chap:Geom::eq:Sobolev_t}
\normlc{\theta}\le C_\theta \normld{\nabla\theta}, \quad\forall \theta\in\Theta.
\end{equation}

Moreover, we consider the tensor space $X^o=Y^o\times \Theta^o\times M^o$, with the following associated norm:
\begin{equation}\label{chap:Geom::eq:Xnorm}
\normX{U}^2=\normld{\nabla\uk}^2+\normld{\nabla\theta}^2+\normld{p}^2, \quad\forall U=(\uk,\theta,p)\in X.
\end{equation}

The variational formulation of problem (\ref{chap:Geom::pb:Cont_Bouss}), over the parameter-dependent original domain is

\begin{equation}\label{chap:Geom::pb:FV}
\left\{\begin{array}{l}
\mbox{Find } (\uk_{o},\theta_{o}^u,p_{o}^u)=(\uk_{o}(\muk),\theta^u_{o}(\muk),p^u_{o}(\muk))\in X^o\mbox{ such that}\vspace{0.3cm}\\

\begin{array}{ll}
\tilde{a}_u(\uk_{o},\vk_{o};\muk)+\tilde{b}(\vk_{o},p^u_{o};\muk)+ \tilde{a}_{Su}'(\uk_{o};\uk_{o},\vk_{o};\muk) \\
+\tilde{c}_u(\uk_{o},\uk_{o},\vk_{o};\muk)+\tilde{f}(\theta_{o}^u,\vk_{o};\muk)=\tilde{F}(\vk_{o};\muk)&\quad\forall\vk_{o}\in Y^o,\vspace{0.15cm}\\

\tilde{b}(\uk_{o},p^v_{o};\muk)=0&\quad\forall p^v_{o}\in M^o,\vspace{0.15cm}\\ 

\tilde{a}_{\theta}(\theta^u_{o},\theta^v_{o};\muk)+\tilde{c}_{\theta}(\uk_{o},\theta_{o}^u,\theta^v_{o};\muk)\\
+\tilde{a}_{S\theta,n}'(\uk_{o};\theta^u_{o},\theta^v_{o};\muk) =\tilde{Q}(\theta^v_{o};\muk)&\quad\forall \theta^v_{o}\in \Theta^o.
\end{array}\end{array}\right.
\end{equation}

Here, the bilinear forms $\tilde{a}_u(\cdot,\cdot;\muk)$, $\tilde{a}_{\theta}(\cdot,\cdot;\muk)$, $\tilde{b}(\cdot,\cdot;\muk)$ and $\tilde{f}(\cdot,\cdot;\muk)$ are defined by 

\begin{equation}\label{eq:bilinear_forms}
\begin{array}{l}
\tilde{a}_u(\uk_o,\vk_o;\muk)=Pr\intOo{\nabla_o\uk_o:\nabla_o\vk_o}, \quad \\
\tilde{a}_\theta(\theta_o^u,\theta_o^v;\muk)=\intOo{\nabla_o\theta_o^u\cdot\nabla_o\theta_o^v},\vspace{0.2cm}\\
\tilde{f}(\theta^u_o,\vk_o;\muk)=-Pr\,\muf\intOo{\theta^u_o\,\vk_{do}},\quad \\
\tilde{b}(\uk_o,p^v_o;\muk)=-\intOo{(\nabla_o\cdot\uk_o)p_o^v};
\end{array}
\end{equation}

the trilinear forms $\tilde{c}_u(\cdot,\cdot,\cdot;\muk)$ and $\tilde{c}_{\theta}(\cdot,\cdot,\cdot;\muk)$ are defined by

\begin{equation}
\begin{array}{l}
\tilde{c}_u(\zk_o,\uk_o,\vk_o;\muk)=\intOo{(\zk_o\cdot\nabla_o\uk_o)\vk_o},\\
\tilde{c}_\theta(\uk_o,\theta^u_o,\theta^v_o;\muk)=\intOo{(\uk_o\cdot\nabla_o\theta^u_o)\theta^v_o};
\end{array}
\end{equation}
and the non-linear Smagorinsky term for eddy viscosity, $\tilde{a}_{Su}'(\cdot;\cdot,\cdot;\muk)$, is given by

\begin{equation}\label{eq:trilinear_forms}
\tilde{a}_{Su}'(\zk_o;\uk_o,\vk_o;\muk)=\intOo{\nu_T(\zk_o)\nabla\uk_o:\nabla\vk_o},
\end{equation}
with
\[
\nu_T(\uk)=(C_Sh_K)^2|\nabla\uk_{|_K}|.
\]
%and $\pih$ an uniformly stable in $H^1$-norm interpolation operator. 

For the thermal eddy diffusivity term $\tilde{a}'_{S\theta,n}$, let us first introduce a mollifier $\phi\in C_c^\infty(\R)$, with supp($\phi$)$\subset B(0,1), \phi\ge0, \|\phi\|_{0,1,\R}>0$, and $\phi$ is even, i.e., $\phi(-x)=\phi(x)$%, \black{and such that $\phi'_n\ra0$ when $n\ra\infty$}
. Let us consider the mollifier sequence $\{\phi_n(x)\}_{n\ge1}$, with $\phi_n\in C_c^\infty(\R)$, supp($\phi_n$)$\subset B(0,1/n),\phi_n\ge0, \|\phi_n\|_{0,1,\R}=1$, defined by
\[
\phi_n(x)=\dfrac{n}{\|\phi\|_{0,1,\R}}\phi(n\,x).
\]
%The following result for mollifiers can be found in \cite{Brezis}:
%\begin{proposition}\label{prop:mollifier}
%Assume $f\in L^p(\R^d)$ with $1\le p<\infty$. Then $(\phi_n*f)\xrightarrow[n\ra\infty]{}f$ in $L^p(\R^n)$.
%\end{proposition}

Thus, the VMS-Smagorinsky eddy diffusivity term is defined as 
\begin{equation}\label{eq:eddy_conduct_reg}
\tilde{a}_{S\theta, n}'(\uk_o;\theta^u_o,\theta^v_o;\muk)
=\intOo{\nu_{T,n}(\uk_o)\nabla_o\theta^{u}_o\cdot\nabla_o\theta^{v}_o},
\end{equation}
with 
\[
\nu_{T,n}(\uk)=(C_Sh_K)^2(\phi_n*|\nabla\uk_{|_K}|),
\]
where $*$ denotes the convolution. \red{The choice of the eddy viscosity is the one suggested in \cite{Smago}. For thermal eddy diffusivity, we have chosen a modified form, considering a mollifier for the thermal eddy diffusivity in \cite{Smago}. The choice of this thermal eddy diffusivity with the mollifier, assure the Lipschitz continuity of the eddy diffusivity operator.}

Thanks to mollifier properties (see \cite{Brezis_en} for details), it holds that $\tilde{a}_{S\theta,n}'$ converges uniformly to $\tilde{a}_{S\theta}'$, with
\[
\tilde{a}_{S\theta}'(\uk_o;\theta_o^u,\theta_o^v;\muk)=\intOo{\nu_{T}(\uk_o)\nabla\theta^{u}_o\cdot\nabla\theta^{v}_o}.
\]

\begin{Remark}
The mollified eddy diffusivity $\tilde{a}_{S\theta, n}'$ is considered for the \\well-possedness analysis and the development of \textit{a posteriori} error estimator. In practice, we consider the eddy diffusivity term in the Boussinesq-Smagorinsky model as $\tilde{a}_{S\theta}'$. 
\end{Remark}

Finally, the linear forms $\tilde{F}$ and $\tilde{Q}$ are given by
\begin{equation}
\tilde{F}(\vk_o;\muk)=\left<\fk,\vk_o\right>, \quad
\tilde{Q}(\theta^v_o;\muk)=\left<Q,\theta^v_o\right>,
\end{equation}
where here $\left<\cdot,\cdot\right>$ stands either the duality paring between $Y$ and $Y'$, and between $\Theta$ and $\Theta'$; being $Y'$ and $\Theta'$ the dual spaces of $Y$ and $\Theta$, respectively.

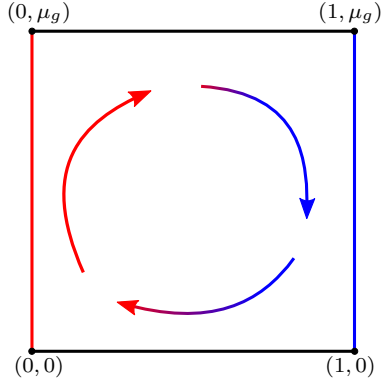
\begin{figure}[h!]
%\pgfdeclareimage[width=0.5\linewidth]{Domain}{Cap_Geom/Cavity_natconv.pdf}
\pgfdeclareimage[width=0.5\linewidth]{Domain}{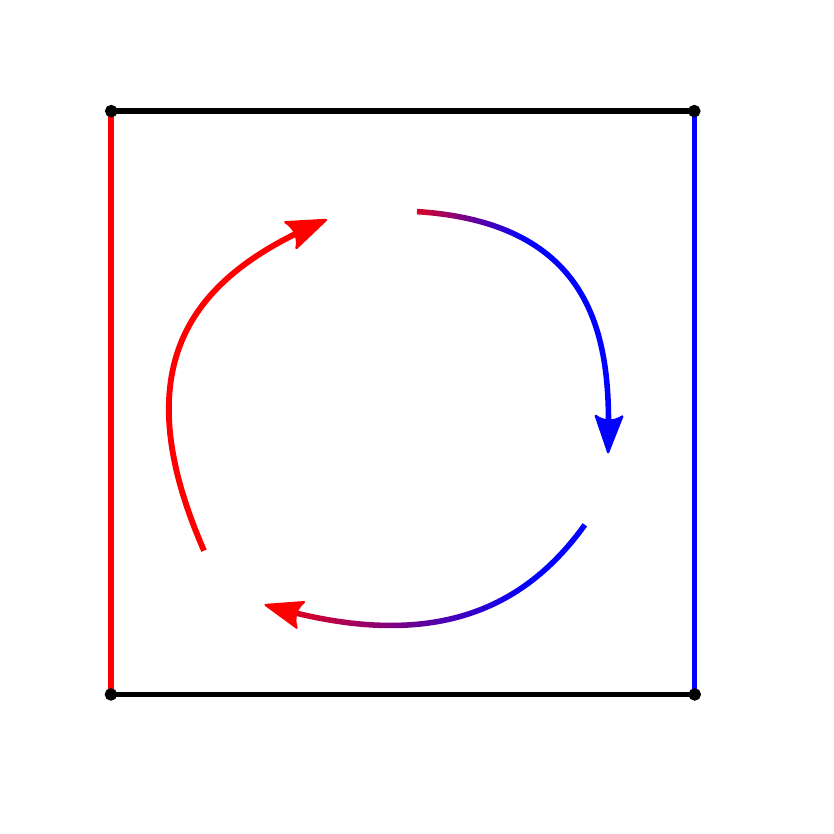}
\centering
  \begin{tikzpicture}[font=\footnotesize]
    \pgftext[at=\pgfpoint{0mm}{0mm},left,base]{\pgfuseimage{Domain}}
    \draw (0.9, 8mm) node {$(0,0)$};
    \draw (5, 8mm) node {$(1,0)$};
    \draw (0.9, 55mm) node {$(0,\mu_g)$};
    \draw (5, 55mm) node {$(1,\mu_g)$};
    %\draw (3, 55mm) node {$\partial_n\theta=0 \quad \uk=0$};
    %\draw (3, 8mm) node {$\partial_n\theta=0 \quad \uk=0$};
    %\draw (0.3, 32mm) node {$\theta=\theta_D$};  
 	%\draw (0.3, 28mm) node {$\uk=0$};
 	%\draw (5.6, 32mm) node {$\theta=0$};  
 	%\draw (5.6, 28mm) node {$\uk=0$};
  \end{tikzpicture}
\caption{Original domain $\Omega_o(\mug)$, with the geometrical parameter considered.}
\label{chap:Geom::fig:problem}
\end{figure}

We present a cavity domain, whose height is varied through a geometrical parameter. This geometrical parameter, that we denote by $\mug$, varies the aspect ratio of the cavity. In Fig. \ref{chap:Geom::fig:problem} we show the original cavity domain considered, with the geometrical parameter considered.

To be able to store parameter independent matrices in the offline phase \black{(see section \ref{chap:Geom::sec:RB})} of the RB method, we need to compute all the integrals in a reference domain trough a transformation of the original domain. Thus, we set $\mug^{ref}=1$, and we define the reference domain $\Omega_r=\Omega_o(\mug^{ref})$. The parameter-dependent original domain can be recovered by a transformation map, $\Psi:\Omega_r\times\cD\ra\R^2$, defined as
\begin{equation}\label{chap:Geom::eq:transformation}
\Psi((x,y);\mug)=\left(\begin{array}{cc}
1 & 0\\
0 & \mug\end{array}\right)\left(\begin{array}{c}x\\y\end{array}\right), \quad\forall (x,y)\in \Omega_r.
\end{equation}

As this map is linear, its Jacobian matrix and its determinant are given by
\begin{equation}\label{chap:Geom::eq:Jacobian}
\textbf{J}((x,y);\mug)=\left(\begin{array}{cc}
1 & 0\\
0 & \mug\end{array}\right), \mbox{ and } |\textbf{J}((x,y);\mug)|=\mug.
\end{equation}

Let $\{\Th\}_{h>0}$ a \black{regular family of triangulations on the reference domain}. We denote $Y_h=(H_0^1(\Omega_r)\cap V_h^l(\Omega_r))^d$, $M_h=L_0^2(\Omega_r)\cap V_h^m(\Omega_r)$, $\Theta_h=H_0^1(\Omega_r)\cap V_h^q(\Omega_r)$, with $l,m,q\in\mathbb{N}$, respectively the discrete velocity, pressure and temperature spaces on the reference domain. \black{Although the analysis can be performed for a regular mesh, for simplicity in the analysis we suppose that the mesh we are considering is uniform.} Denoting $X_h=Y_h\times\Theta_h\times M_h$, we rewrite problem (\ref{chap:Geom::pb:FV}) with respect to the reference domain, applying the change of variables of the transformation map $T$, as 
\begin{equation}\label{chap:Geom::pb:FVref}
\left\{\begin{array}{l}
\mbox{Find } (\uk_{h},\theta_{h}^u,p_{h}^u)=(\uk_{h}(\muk),\theta^u_{h}(\muk),p^u_{h}(\muk))\in X_h\mbox{ such that}\vspace{0.3cm}\\

\begin{array}{ll}
a_{u,x}(\uk_h,\vk_h;\muk)+a_{u,y}(\uk_h,\vk_h;\muk)+b_{x}(\vk_h,p^u_h;\muk)\\
+b_{y}(\vk_h,p^u_h;\muk)+ a_{Su,x}'(\uk_h;\uk_h,\vk_h;\muk) \\
+a_{Su,y}'(\uk_h;\uk_h,\vk_h;\muk)+c_{u,x}(\uk_h,\uk_h,\vk_h;\muk)\\
+c_{u,y}(\uk_h,\uk_h,\vk_h;\muk)+f(\theta_{h}^u,\vk_{h};\muk)=F(\vk_{h};\muk)&\quad\forall\vk_{h}\in Y_h,\vspace{0.15cm}\\

b_{x}(\uk_{h},p^v_{h};\muk)+b_{y}(\uk_{h},p^v_{h};\muk)=0&\quad\forall p^v_{h}\in M_h,\vspace{0.15cm}\\ 

a_{\theta,x}(\theta^u_{h},\theta^v_{h};\muk)+a_{\theta,y}(\theta^u_{h},\theta^v_{h};\muk)+c_{\theta,x}(\uk_{h},\theta_{h}^u,\theta^v_{h};\muk)\\
+c_{\theta,y}(\uk_{h},\theta_{h}^u,\theta^v_{h};\muk)+a_{S\theta,nx}'(\uk_{h};\theta^u_{h},\theta^v_{h};\muk)\\
+a_{S\theta,ny}'(\uk_{h};\theta^u_{h},\theta^v_{h};\muk) =Q(\theta^v_{h};\muk)&\quad\forall \theta^v_{h}\in \Theta_h,
\end{array}\end{array}\right.
\end{equation}
where the subscripts $x$ and $y$ \black{denote the addend} of the corresponding operator, \black{relative to the partial derivative} with respect to $x$ or $y$, respectively. \black{For both the eddy viscosity and eddy diffusivity terms, we consider a VMS small-small setting approach (\textit{cf.} \cite{LES-VMS}). For that, \black{we consider a uniformly} $H^1$-norm interpolator operator from $Y_h$ to $\overline{Y}_h=(H_0^1(\Omega_r)\cap V_h^{l-1}(\Omega_r))^d$,  denoted by $\Pi_h$. Thus, we assume that $\Pi_h$ satisfies that there exists a constant $C_f>0$ independent of $h$ such that
\begin{equation}\label{eq:pih_ineq}
\normh{\pih\uk_h}\le C_f\normh{\uk_h} \quad \uk_h\in Y_h,
\end{equation}
where we denote $\pih=Id-\Pi_h$. See \cite{TomasSmago,EDelgado_PhD} for more details.}

The operators in (\ref{chap:Geom::pb:FVref}) have the following form:

\[
a_{u,x}(\uk_h,\vk_h;\muk)=Pr\,\mug\intOr{(\partial_xu_1\partial_xv_1+\partial_xu_2\partial_xv_2)},
\]
\[
a_{u,y}(\uk_h,\vk_h;\muk)=\dfrac{Pr}{\mug}\intOr{(\partial_yu_1\partial_yv_1+\partial_yu_2\partial_yv_2)},
\]
\[
b_{x}(\uk_h,p^u_h;\muk)=-\mug\intOr{p^u_h\partial_xu_1}, \quad 
b_{y}(\uk_h,p^u_h;\muk)=-\intOr{p^u_h\partial_yu_2},
\]
\[
f(\theta_{h}^u,\vk_{h};\muk)=Pr\,\muf\,\mug\intOr{\theta^u_hv_2},
\]
\[
a_{\theta,x}(\theta^u_{h},\theta^v_{h};\muk)=\mug\intOr{\partial_x\theta_h^u\partial_x\theta_h^v},\;\;
a_{\theta,y}(\theta^u_{h},\theta^v_{h};\muk)=\dfrac{1}{\mug}\intOr{\partial_y\theta_h^u\partial_y\theta_h^v},
\]
\[
c_{u,x}(\wk_h,\uk_h,\vk_h;\muk)=\mug\intOr{\big[(w_1\partial_xu_1)v_1+(w_1\partial_xu_2)v_2\big]},
\]
\[
c_{u,y}(\wk_h,\uk_h,\vk_h;\muk)=\intOr{\big[(w_2\partial_yu_1)v_1+(w_2\partial_yu_2)v_2\big]},
\]
\[
c_{\theta,x}(\uk_{h},\theta_{h}^u,\theta^v_{h};\muk)=\mug\intOr{(u_1\partial_x\theta^u_h)\theta^v_h},
\]
\[
c_{\theta,y}(\uk_{h},\theta_{h}^u,\theta^v_{h};\muk)=\intOr{(u_2\partial_y\theta^u_h)\theta_h^v},
\]
\[
a_{Su,x}'(\wk_h;\uk_h,\vk_h;\muk)=\mug\hspace{-0.15cm}\int_{\Omega_r}{\hspace{-0.2cm}\nu_T(\pih\wk)\big[\partial_x(\pih u_1)\partial_x(\pih v_1)}\]\[
{+\partial_x(\pih u_2)\partial_x(\pih v_2)\dOr\big]},
\]
\[
a_{Su,y}'(\wk_h;\uk_h,\vk_h;\muk)=\dfrac{1}{\mug}\hspace{-0.1cm}\int_{\Omega_r}{\hspace{-0.2cm}\nu_T(\pih\wk)\big[\partial_y(\pih u_1)\partial_y(\pih v_1)}\]\[
{+\partial_y(\pih u_2)\partial_y(\pih v_2)\dOr\big]},
\]
\[
a_{S\theta,nx}'(\uk_{h};\theta^u_{h},\theta^v_{h};\muk)=\dfrac{\mug}{Pr}\intOr{\nu_{T,n}(\pih\uk;\mug)\partial_x(\pih \theta^u)\partial_x(\pih \theta^v)},
\]
\[
a_{S\theta,ny}'(\uk_{h};\theta^u_{h},\theta^v_{h};\muk)=\dfrac{1}{Pr\,\mug}\intOr{\nu_{T,n}(\pih\uk;\mug)\partial_y(\pih \theta^u)\partial_y(\pih \theta^v)},
\]

These integrals are derived applying the well-known change of variable formula (see e.g. \cite{LibroBR}). With this geometrical parametrization, the eddy viscosity $\nu_T(\cdot)$ (analogously $\nu_{T,n}(\cdot)$) also depends on the geometrical parameter, and is defined as
\begin{equation}\label{chap:Geom::eq:eddy_visco}
\nu_T(\uk;\mug)=C_S^2\frac{\mug^2+1}{N_h^2}\sqrt{(\partial_x u_1)^2+\dfrac{1}{\mug^2}(\partial_yu_1)^2+(\partial_xu_2)^2+\dfrac{1}{\mug^2}(\partial_yu_2)^2}.
\end{equation}

Here we are supposing that we consider an uniform mesh in the reference domain $\Omega_r$, with $N_h$ partitions on each side. Since the mesh size, $h_K$, in the VMS-Smagorinsky eddy viscosity and eddy diffusivity appears in (\ref{chap:Geom::eq:eddy_visco}) in terms of the parameter-dependent original domain, we map it to the reference domain, by applying the change of variable map $\Psi$ defined in (\ref{chap:Geom::eq:transformation}).%, obtaining that

%\[h_K^2=\dfrac{\mug^2+1}{N_h^2},\quad \forall K\in\Th.\]

%Note that problem (\ref{chap:Bouss::pb:FVX}) is a particular case of problem (\ref{chap:Geom::pb:FVref}), taking $\mug=1$.

\section{Reduced Basis formulation}\label{chap:Geom::sec:RB}
In this section we present the RB problem derived from the discrete problem presented in section \ref{chap:Geom::sec:FE}. We construct the low-dimensional spaces for the RB problem with the Greedy algorithm. Both the pressure and the temperature reduced basis spaces are defined with the corresponding snapshots computed solving the FE problem (\ref{chap:Geom::pb:FVref}).

The RB velocity space is constructed with the velocity snapshot of the FE velocity solution, and the inner pressure \textit{supremizer} \black{(\textit{cf.} \cite{Rozza2006,Stokes4})}, $T_p^{\muk}: M_h\ra Y_h$, defined for this problem as
\begin{equation}
\intOr{\nabla T_p^{\muk} q_h:\nabla\vk_h}=-\mug\intOr{q_h\,\partial_xv_1}-\intOr{q_h\,\partial_yv_2}, \quad \forall \vk_h\in Y_h.
\end{equation}

Thus, the reduced basis spaces are given by
\begin{equation}\label{chap:Geom::YN}
Y_{N}=\mbox{span}\{\zetav_{2k-1}:=\uk_h(\muk^k),~ \zetav_{2k}:=T_p^{\muk}\xip_k, ~ ~ k=1,\dots,N\},
\end{equation}
\begin{equation}\label{chap:Geom::M_N}
M_{N}=\mbox{span}\{\xip_k:=p_h^u(\muk^k),~ ~ k=1,\dots,N\},
\end{equation}
\begin{equation}\label{chap:Geom::M_N}
\Theta_{N}=\mbox{span}\{\phit_k:=\theta^u_h(\muk^k),~ ~ k=1,\dots,N\}.
\end{equation}

Denoting $X_N=Y_N\times\Theta_N\times M_N$, the RB problem is
\begin{equation}\label{chap:Geom::pb:RBPref}
\left\{\begin{array}{l}
\mbox{Find } (\uk_N(\muk),\theta_N^u(\muk),p_N^u(\muk)\in X_N\mbox{ such that}\vspace{0.3cm}\\

\begin{array}{ll}
a_{u,x}(\uk_N,\vk_N;\muk)+a_{u,y}(\uk_N,\vk_N;\muk)+b_{x}(\vk_N,p^u_N;\muk)\\
+b_{y}(\vk_N,p^u_N;\muk)+ a_{Su,x}'(\uk_N;\uk_N,\vk_N;\muk) \\
+a_{Su,y}'(\uk_N;\uk_N,\vk_N;\muk)+c_{u,x}(\uk_N,\uk_N,\vk_N;\muk)\\
+c_{u,y}(\uk_N,\uk_N,\vk_N;\muk)+f(\theta_{N}^u,\vk_{N};\muk)=F(\vk_{N};\muk)&\quad\forall\vk_{N}\in Y_N,\vspace{0.15cm}\\

b_{x}(\uk_{N},p^v_{N};\muk)+b_{y}(\uk_{N},p^v_{N};\muk)=0&\quad\forall p^v_{N}\in M_N,\vspace{0.15cm}\\ 

a_{\theta,x}(\theta^u_{N},\theta^v_{N};\muk)+a_{\theta,y}(\theta^u_{N},\theta^v_{N};\muk)+c_{\theta,x}(\uk_{N},\theta_{N}^u,\theta^v_{N};\muk)\\
+c_{\theta,y}(\uk_{N},\theta_{N}^u,\theta^v_{N};\muk)+a_{S\theta,nx}'(\uk_{N};\theta^u_{N},\theta^v_{N};\muk)\\
+a_{S\theta,ny}'(\uk_{N};\theta^u_{N},\theta^v_{N};\muk) =Q(\theta^v_{N};\muk)&\quad\forall \theta^v_{N}\in \Theta_N.
\end{array}\end{array}\right.
\end{equation}

The eddy viscosity $\nu_T(\uk;\mug)$ must be \black{tensorized} in problem (\ref{chap:Geom::pb:RBPref}), for the efficient solve in the online phase. For this purpose, we consider the use of EIM (\textit{cf. }\cite{EIM2,EIM1}). The eddy viscosity and eddy diffusivity terms are approximated as
\[
a_{Su,x}'(\wk_h;\uk_h,\vk_h;\muk)\approx\hat{a}_{Su,x}'(\uk_N,\vk_N;\muk),
\]
\[
a_{Su,y}'(\wk_h;\uk_h,\vk_h;\muk)\approx\hat{a}_{Su,y}'(\uk_N,\vk_N;\muk),
\]
\[
a_{S\theta,nx}'(\uk_{h};\theta^u_{h},\theta^v_{h};\muk)\approx\hat{a}_{S\theta,x}'(\theta^u_N,\theta^v_N;\muk),
\]
\[
a_{S\theta,ny}'(\uk_{h};\theta^u_{h},\theta^v_{h};\muk)\approx\hat{a}_{S\theta,y}'(\theta^u_N,\theta^v_{N};\muk),
\]
with,
\[
\hat{a}_{Su,x}'(\uk_N,\vk_N;\muk)=\sum_{k=1}^M\sigma_k(\muk)s_{u,x}(q_k,\uk_N,\vk_N),
\]
\[
\hat{a}_{Su,y}'(\uk_N,\vk_N;\muk)=\sum_{k=1}^M\sigma_k(\muk)s_{u,y}(q_k,\uk_N,\vk_N),
\]
\[
\hat{a}_{S\theta,x}'(\theta^u_N,\theta^v_N;\muk)=\dsum_{k=1}^{M}\sigma_k(\muk)s_{\theta,x}(q_k,\theta^u_N,\theta^v_N),
\]
\[
\hat{a}_{S\theta,y}'(\theta^u_N,\theta^v_N;\muk)=\dsum_{k=1}^{M}\sigma_k(\muk)s_{\theta,y}(q_k,\theta^u_N,\theta^v_N),
\]
and,
\[
s_{u,x}(q_k,\uk_N,\vk_N)=\mug\intOr{q_k \,\big[\partial_x(\pih u_1)\partial_x(\pih v_1)
+\partial_x(\pih u_2)\partial_x(\pih v_2)\big]},
\]
\[
s_{u,y}(q_k,\uk_N,\vk_N)=\frac{1}{\mug}\intOr{q_k \,\big[\partial_y(\pih u_1)\partial_y(\pih v_1)
+\partial_y(\pih u_2)\partial_y(\pih v_2)\big]},
\]
\[
s_{\theta,x}(q_k,\theta^u_N,\theta^v_N)=\dfrac{\mug}{Pr}\intOr{q_k\,\partial_x(\pih \theta^u)\partial_x(\pih \theta^v)},
\]
\[
s_{\theta,y}(q_k,\theta^u_N,\theta^v_N)=\dfrac{1}{Pr\,\mug}\intOr{q_k\,\partial_y(\pih \theta^u)\partial_y(\pih \theta^v)},
\]
\black{where $\sigma_k(\muk)$ and $q_k$ are computed by the EIM algorithm (see \cite{EIM1} for further details).} 

The parameter independent matrices and tensors to store during the offline phase in order to efficiently solve  problem (\ref{chap:Geom::pb:RBPref}), are given in this case by
\[
(\mathbb{A}_N^{u,x})_{ij}=a_{u,x}(\zeta^{\vk}_j,\zeta^{\vk}_i),\enspace (\mathbb{A}_N^{u,y})_{ij}=a_{u,y}(\zeta^{\vk}_j,\zeta^{\vk}_i),\quad i,j=1,\dots,2N,
\]
\[
(\mathbb{A}_N^{\theta,x})_{lm}=a_{\theta,x}(\varphi^{\theta}_m,\phit_l),\enspace (\mathbb{A}_N^{\theta,y})_{lm}=a_{\theta,y}(\varphi^{\theta}_m,\phit_l),\quad  l,m=1,\dots,N,
\]
\[
(\mathbb{F}_N)_{li}=f(\phit_l,\zeta^{\vk}_i),
\quad i=1,\dots,2N,\;l=1,\dots,N,
\]
\[
(\mathbb{B}_N^{x})_{li}=b_{x}(\zeta^{\vk}_i,\xi^p_l),\enspace(\mathbb{B}_N^{(y)})_{li}=b_{y}(\zeta^{\vk}_i,\xi^p_l),\quad i=1,\dots,2N,\;l=1,\dots,N,
\]
\[
(\mathbb{C}_N^{u,x}(\zeta^{\vk}_s))_{ij}=c_{u,x}(\zeta^{\vk}_s,\zeta^{\vk}_j,\zeta^{\vk}_i), 
\quad i,j,s=1,\dots, 2N, 
\]
\[
(\mathbb{C}_N^{u,y}(\zeta^{\vk}_s))_{ij}=c_{u,y}(\zeta^{\vk}_s,\zeta^{\vk}_j,\zeta^{\vk}_i), 
\quad i,j,s=1,\dots, 2N, 
\]
\[
(\mathbb{C}_N^{\theta,x}(\zeta^{\vk}_s))_{lm}=c_{\theta,x}(\zeta^{\vk}_s,\phit_m,\phit_l), 
\quad l,m=1,\dots,N,\;s=1,\dots, 2N,  
\]
\[
(\mathbb{C}_N^{\theta,y}(\zeta^{\vk}_s))_{lm}=c_{\theta,y}(\zeta^{\vk}_s,\phit_m,\phit_l), 
\quad l,m=1,\dots,N,\;s=1,\dots, 2N,  
\]
\[
(\mathbb{S}^{u,x}_N(q_s))_{ij}=s_{u,x}(q_s,\zeta^{\vk}_j,\zeta^{\vk}_i),
\quad i,j=1,\dots,2N, s=1,\dots, M,
\]
\[
(\mathbb{S}^{u,y}_N(q_s))_{ij}=s_{u,y}(q_s,\zeta^{\vk}_j,\zeta^{\vk}_i),
\quad i,j=1,\dots,2N, s=1,\dots, M,
\]
\[
(\mathbb{S}^{\theta,x}_N(q_s))_{lm}=s_{\theta,x}(q_s,\phit_m,\phit_l),
\quad l,m=1,\dots,N, s=1,\dots, M, 
\]
\[
(\mathbb{S}^{\theta,y}_N(q_s))_{lm}=s_{\theta,y}(q_s,\phit_m,\phit_l),
\quad l,m=1,\dots,N, s=1,\dots, M.  
\]

Here we are representing the reduced basis velocity, temperature and pressure solutions as a linear combination of the velocity, temperature and pressure snapshots, respectively, of the reduced spaces, i.e.,
\[
\uk_N(\muk)=\sum_{j=1}^{2N}u_j^N(\muk)\zeta^{\vk}_j,
\quad \theta_N(\muk)=\sum_{j=1}^N\theta_j^N(\muk)\phit_j,
\quad p_N(\muk)=\sum_{j=1}^Np_j^N(\muk)\xi^p_j.
\]

\section{\textit{A posteriori} error estimator}\label{chap:Geom::sec:a_post}
In order to develop the \textit{a posteriori} error estimator for the Greedy algorithm, we rewrite problem (\ref{chap:Geom::pb:FVref}) in a more compact form as
\begin{equation}\label{chap:Geom::pb:FVX}\left\{\begin{array}{l}
\mbox{Find }U_h(\muk)=(\uk_h,\theta_h^u, p_h^u)\in X_h \mbox{ such that}\vspace{0.2cm}\\
A(U_h(\muk),V_h;\muk)=F(V_h;\muk)\qquad\forall V_h\in X_h,\end{array}\right.
\end{equation}

%\begin{remark}
%For the numerical analysis in the development of the \textit{a posteriori } error bound estimator, we consider the regularized, by a mollifier, eddy diffusivity defined in (\ref{chap:Bouss::eq:eddy_conduct_reg}).
%\end{remark}

The \textit{a posteriori} error estimator is based upon the BRR theory (\textit{cf. }\cite{BRR}). For this purpose we define the Gateaux derivative of $A(\cdot,\cdot;\muk)$ with respect to the first variable, in the direction $Z\in X$, denoted by $\partial_1 A(U,V;\muk)(Z)$. For this problem, \black{denoting $Z=(\zk,\theta^z,p^z)$}, the derivative is defined by:

\[
\partial_1 A(U,V;\muk)(Z)=a_{u,x}(\zk,\vk;\muk)+a_{u,y}(\zk,\vk;\muk)+b_{x}(\vk,p^z;\muk)
+b_{y}(\vk,p^z;\muk)
\]
\[
+f(\theta^z,\vk;\muk)-b_{x}(\zk,p^v;\muk)-b_{y}(\zk,p^v;\muk)+a_{\theta,x}(\theta^z,\theta^v;\muk)
+a_{\theta,y}(\theta^z,\theta^v;\muk)
\]
\[
+c_{u,x}(\zk,\uk,\vk;\muk)+c_{u,x}(\uk,\zk,\vk;\muk)
+c_{u,y}(\zk,\uk,\vk;\muk)+c_{u,y}(\uk,\zk,\vk;\muk)
\]
\[
+c_{\theta,x}(\zk,\theta^u,\theta^v;\muk)+c_{\theta,x}(\uk,\theta^z,\theta^v;\muk)
+c_{\theta,y}(\zk,\theta^u,\theta^v;\muk)+c_{\theta,y}(\uk,\theta^z,\theta^v;\muk)
\]
\[
+a_{Su,x}'(\uk;\zk,\vk;\muk) +a_{Su,y}'(\uk;\zk,\vk;\muk)
+a_{S\theta,nx}'(\uk;\theta^z,\theta^v;\muk)
+a_{S\theta,ny}'(\uk;\theta^z,\theta^v;\muk)
\]
\[
+\mug\intOr{\partial_1\nu_T(\pih \uk)(\pih \zk)[\dx(\pih u_1)\dx(\pih v_1)+\dx(\pih u_2)\dx(\pih v_2)]}
\]
\[
+\dfrac{1}{\mug}\intOr{\partial_1\nu_T(\pih \uk)(\pih \zk)[\dy(\pih u_1)\dy(\pih v_1)+\dy(\pih u_2)\dy(\pih v_2)]}
\]
\[
+\dfrac{\mug}{Pr}\intOr{\partial_1\nu_{T,n}(\pih \uk)(\pih \zk)\;\dx(\pih \theta^u)\dx(\pih \theta^v)}
\]
\[
+\dfrac{1}{\mug\,Pr}\intOr{\partial_1\nu_{T,n}(\pih \uk)(\pih \zk)\;\dy(\pih \theta^u)\dy(\pih \theta^v)},
\]
with
\[
\partial_1\nu_T(\uk)(\zk)=C_S^2\frac{\mug^2+1}{N_h^2}\dfrac{\partial_x u_1\dx z_1+\dfrac{1}{\mug^2}\partial_yu_1\dy z_1+\partial_xu_2\dx z_2+\dfrac{1}{\mug^2}\partial_yu_2\dy z_2}{|\nabla(\Psi^{-1}\uk)|},
%\sqrt{(\partial_x u_1)^2+\dfrac{1}{\mug^2}(\partial_yu_1)^2+(\partial_xu_2)^2+\dfrac{1}{\mug^2}(\partial_yu_2)^2}
\]
and
\[
\partial_1\nu_{T,n}(\uk)(\zk)=C_S^2\frac{\mug^2+1}{N_h^2}\left[\phi_n'*|\nabla \Psi^{-1}(\uk)|\right]:\left[\nabla \Psi^{-1}(\zk)\right].
%\left[\phi_n*\sqrt{(\partial_x u_1)^2+\dfrac{1}{\mug^2}(\partial_yu_1)^2+(\partial_xu_2)^2+\dfrac{1}{\mug^2}(\partial_yu_2)^2}\right]
\]
%\[
%\cdot(\dx z_1+\frac{1}{\mug^2}\dy z_1+\dx z_2+\frac{1}{\mug^2}\dy z_2).
%\]

The Gateaux derivative satisfies the following continuity and inf-sup conditions:
\begin{equation}\label{chap:Geom::ec:cont}
\infty>\gamma_0\ge\gamma_h(\muk)\equiv\sup_{Z_h\in X_h}\sup_{V_h\in X_h}\dfrac{\partial_1A(U_h(\mu),V_h;\muk)(Z_h)}{\|Z_h\|_X\|V_h\|_X}.
\end{equation}
\begin{equation}\label{chap:Geom::ec:infsup}
0<\beta_0<\beta_{h}(\muk)\equiv\inf_{Z_h\in X_h}\sup_{V_h\in X_h}\dfrac{\partial_1A(U_h(\mu),V_h;\muk)(Z_h)}{\|Z_h\|_X\|V_h\|_X}.
\end{equation}

The existence of $\gamma_0\in\R$ and $\beta_0>0$ satisfying (\ref{chap:Geom::ec:cont}) and (\ref{chap:Geom::ec:infsup}), respectively, are given by the following results, \red{whose proofs can be found in \ref{app:prop1} and \ref{app:prop2} respectively}.

\begin{proposition}\label{chap:Geom::prop::cont}
There exists $\gamma_0\in\R$ such that $\forall \mu\in\cD$
\[
|\partial_1A(U_h(\muk),V_h;\muk)(Z_h)|\le \gamma_0\|Z_h\|_X\|V_h\|_X\quad \forall Z_h,V_h\in X_h.
\]
\end{proposition}

\begin{proposition}\label{chap:Geom::prop:infsup}
\black{Let $C(\muk,\normlur{\phi_n'})$ a constant depending on $\muk$ and $\normlur{\phi_n'}$. Then if 
\begin{equation}\label{infsup::vel}
\normld{\nabla\uk_h}\le\dfrac{2Pr\min\{\mug,1/\mug\}-C_PPr\mug\mu}{4C_u^2\min\{\mug,1\}+C(\muk,\normlur{\phi_n'})\normld{\nabla\theta_h^u}}
\end{equation}
and
\begin{equation}\label{infsup::temp}
\normld{\nabla\theta_h^u}\le\dfrac{2\min\{\mug,1/\mug\}-C_PPr\mug\mu}{4C_uC_{\theta}\min\{\mug,1\}+C(\muk,\normlur{\phi_n'})\normld{\nabla\uk_h}},
\end{equation}}
then there exists $\tilde{\beta}(\muk)>0$ such that
\begin{equation}\label{chap:Geom::eq:infsup_prop}
\partial_1A(U_h,V_h;\muk)(V_h)\ge\tilde{\beta}(\muk)(\normld{\nabla\vk_h}^2+\normld{\nabla\theta_h^v}^2)\quad\forall V_h\in\ X_h.
\end{equation}
\end{proposition}

Since operator $b(\vk_h,p_h;\muk)$ satisfies the discrete inf-sup condition
\[
\alpha\normld{p_h}\le \sup_{\vk_h\in Y_h}\dfrac{b(\vk_h,p_h;\muk)}{\normld{\vk_h}},
\]

we can prove that the inf-sup condition (\ref{chap:Geom::eq:infsup_prop}) is satisfied thanks to Prop. \ref{chap:Geom::prop:infsup}.
%Proposition \ref{chap:Geom::prop::cont} and Proposition \ref{chap:Geom::prop:infsup}, can be proved analogously as in Proposition \ref{chap:Bouss::prop::cont} and Proposition \ref{chap:Bouss::prop:infsup}, respectively, taking into account that in the derivative of the operator $A(\cdot,\cdot;\muk)$ of problem (\ref{chap:Geom::pb:FVX}) there is a re-scaling of the terms in the derivative operator  $A(\cdot,\cdot;\mu)$ of problem (\ref{chap:Bouss::pb:FVX}).
\newline

With the following result, we prove that the Gateaux derivative of the Boussinesq-Smagorinsky operator is locally Lipschitz-continuous. \red{The proof can be found in \ref{app:Lemma}.}

\begin{lemma}\label{chap:Geom::lema:LemmaRho}
Let $U_h^1,U_h^2\in X_h$. Then, in a neighborhood of $U_h^1$ and $U_h^2$, there exists a positive constant $\rho_n(\mug)$ such that, $\forall Z_h,V_h \in X_h$,
\begin{equation}\label{chap:Geom::ec:ro}
\left|\partial_1A(U_h^1,V_h;\muk)(Z_h)-\partial_1A(U_h^2,V_h;\muk)(Z_h)\right|\le\rho_n(\mug)\normX{U_h^1-U_h^2}\normX{Z_h}\|V_h\|_X.
\end{equation}
%with
%\[\rho_n(\mug)=2\max\{\mug,1\}C_u+2\max\{\mug,1\}C_{\theta}+4\max\left\{\mug,\frac{1}{\mug}\right\}C_S^2h^{2-d/2}CC_f^3\]
%\[+\max\left\{\mug,\frac{1}{\mug}\right\}\dfrac{C_S^2h^{2-d/2}CC_f^3\normlur{\phi_n'}(\normlinf{\nabla(\pih\uk_h^1)}+\normlinf{\nabla(\pih\theta^{u2})})}{Pr}.\]
\end{lemma}

We define the following continuity and inf-sup constants: 
\begin{equation}\label{chap:Geom::betaN}
0<\beta_{N}(\muk)\equiv\inf_{Z_h\in X_h}\sup_{V_h\in X_h}\dfrac{\partial_1A(U_N(\muk),V_h;\muk)(Z_h)}{\|Z_h\|_X\|V_h\|_X}=\inf_{Z_h\in X_h}\dfrac{\|T_N Z_h\|_X}{\|Z_h\|_X},
\end{equation}
\begin{equation}\label{chap:Geom::gammaN}
\infty>\gamma_{N}(\muk)\equiv\sup_{Z_h\in X_h}\sup_{V_h\in X_h}\dfrac{\partial_1A(U_N(\muk),V_h;\muk)(Z_h)}{\|Z_h\|_X\|V_h\|_X}=\sup_{Z_h\in X_h}\dfrac{\|T_N Z_h\|_{X}}{\|Z_h\|_X},
\end{equation}
where the \textit{supremizer} operator $T_N$ is defined as 
\begin{equation}\label{chap:Smago::TN}
(T_NZ_h,V_h)_X=\partial_1A(U_N(\mu),V_h;\mu)(Z_h)\quad \forall V_h,Z_h\in X_h,
\end{equation}
such that
\begin{equation}
T_NZ_h=\arg\sup_{V_h\in X_h}\dfrac{\partial_1A(U_N(\mu),V_h;\mu)(Z_h)}{\|V_h\|_X}.
\end{equation}

The existence of these constants can be proved in the same way that the existence of the constants (\ref{chap:Geom::ec:cont})-(\ref{chap:Geom::ec:infsup}).  Thus, we can define the \textit{a posteriori} error estimator as
\begin{equation}\label{chap:Geom::delta}
\dnk=\frac{\beta_N(\muk)}{2\rho_n(\mug)}\left[1-\sqrt{1-\taunk}\right],
\end{equation}
where $\taunk$ is given by
\begin{align}
\taunk&=\frac{4\enk\rho_n(\mug)}{\beta_N^2(\muk)},\label{chap:Geom::tau}
\end{align} 
with $\enk$ the dual norm of the residual. The \textit{a posteriori} error estimator is stated by the following result, \red{whose proof can be found in \ref{app:Teorema}}. 

\begin{theorem}\label{chap:Geom::teor:Teorprinc}
Let $\muk\in\cD$, and assume that $\beta_N(\muk)>0$. If problem (\ref{chap:Geom::pb:FVX}) admits a solution $U_h(\muk)$ such that
\[
\normX{U_h(\muk)-U_N(\muk)}\le\frac{\beta_N(\muk)}{\rho_n(\mug)},
\]
then this solution is unique in the ball $B_X\left(U_N(\muk),\dfrac{\beta_N(\muk)}{\rho_n(\mug)}\right)$. 

Moreover, assume that $\taunk\le1$ for all $\muk\in\cD$. Then there exists a unique solution $U_h(\muk)$ of (\ref{chap:Geom::pb:FVX}) such that the error with respect $U_N(\muk)$, solution of (\ref{chap:Geom::pb:RBPref}), is bounded by the \textit{a posteriori} error estimator, i.e.,
\begin{equation}\label{chap:Geom::teorprinc:err}
\normX{U_h(\muk)-U_N(\muk)}\le\dnk,
\end{equation}
with effectivity
\begin{equation}\label{chap:Geom::teorprinc:efec}
\dnk\le\left[\frac{2\gamma_{N}(\muk)}{\beta_N(\muk)}+\taunk\right]\normX{U_h(\muk)-U_N(\muk)}.
\end{equation}
\end{theorem}

\section{Numerical results}\label{chap:Geom::sec:Num}

In this section we present some numerical results for the Boussinesq VMS-Smago\-rinsky RB model. We perform \black{three} different configurations for the parametrical set. The first configuration \red{corresponds} to the consideration only of physical parametrical set, fixing the value of the geometrical parameter $\mug=1$. %In this two first configurations, we start with moderate Rayleigh number range, with $Ra\in[10^3,10^5]$, and then we increase the value of Rayleigh number, considering higher Rayleigh number values ranging in $[10^5,10^6]$. 
\black{Here we consider two different scenarios depending on the Rayleigh number range. In order to get error levels small enough for taking into account the \textit{a posteriori} error estimator, we split the Rayleigh number ranger considered, $\muf\in[10^3,10^6]$, into two ranges, $\muf\in[10^3,10^5]$ and $\muf\in[10^5,10^6]$.}

Then, we suppose that the Rayleigh number is fixed with $\muf=10^5$, and we consider the geometrical parameter ranging in $\mug\in\cD=[0.5,2]$. Finally, we consider both the geometrical parameter and the Rayleigh number. For this test, we consider the Rayleigh number, $\muf$, ranging in $[10^3,10^4]$, and the geometrical parameter, $\mug$, ranging in $\mug\in[0.5,2]$. Thus we are considering that the parameter domain is $\cD=[10^3,10^4]\times[0.5,2]$. For all cases, the Prandtl number considered is $Pr=0.71$, that corresponds to the air Prandtl number.

In all tests, we consider no-slip boundary conditions for velocity. We also consider homogeneous Neumann boundary condition for temperature at the top and bottom of the cavity, and Dirichlet conditions for the vertical walls: $\theta=1$ for the left vertical wall and $\theta=0$ for the right vertical wall. \black{Moreover, we consider $\fk=0$ and $Q=0$ both for the momentum equation and the energy equation, respectively.}

The FE solution is computed through a semi-implicit evolution approach, considering that the steady state solution is reached when the error between two iterates is below $\veps_{FE}=10^{-10}$. \black{The FE solution has been computed} considering $\mathbb{P}2-\mathbb{P}2-\mathbb{P}1$ finite elements for velocity, temperature and pressure, respectively.

\black{With respect to the constants that appear in the \textit{a posteriori } error estimator, we explain in the following the numerical approximation of those constants. For the constant $\beta$, we use the Radial Basis Function (RBF) algorithm in order to compute efficiently $\beta(\mu)$ for all $\mu \in \cal{D}$. We follow the technique suggested in \cite{infsup}. Although it can not be proved that a lower bound of the constant $\beta$ as the Successive Constraint Method (SCM) is provided, the computational time is much lower for the RBF, mainly when more than one parameter is considered, with quite good accuracy. The constant $\rho_n$, and more precisely the Sobolev embedding constants, are built once in the reference domain following the algorithm proposed in \cite{Deparis}. Moreover, we compute the Lipschitz constant $\rho_n(\mu_g)$ without taking into account the term in which the mollifier takes part of it. This term in $\rho_n(\mu_g)$ that comes from the mollifier is multiplied by a factor of $(C_Sh)^2$, that has no relevance in the value of the Lipschitz constant. Thus, we consider accurate the approximation of the Lipschitz constant without considering the mollifier.  }

\subsection{Physical parametrization}\label{chap:Bouss::sec:High_Ra}
\begin{figure}[h!]
\centering
\includegraphics[width=0.3\linewidth]{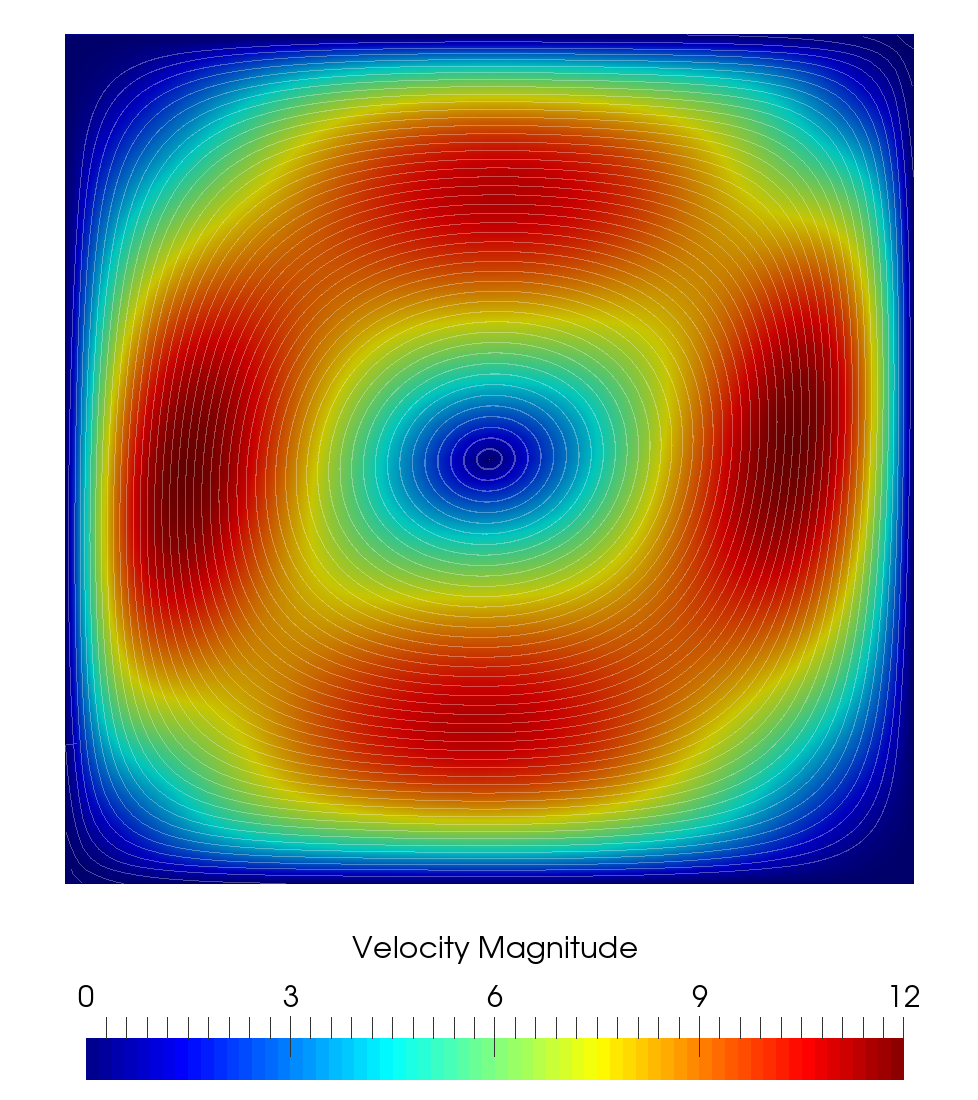}
\includegraphics[width=0.3\linewidth]{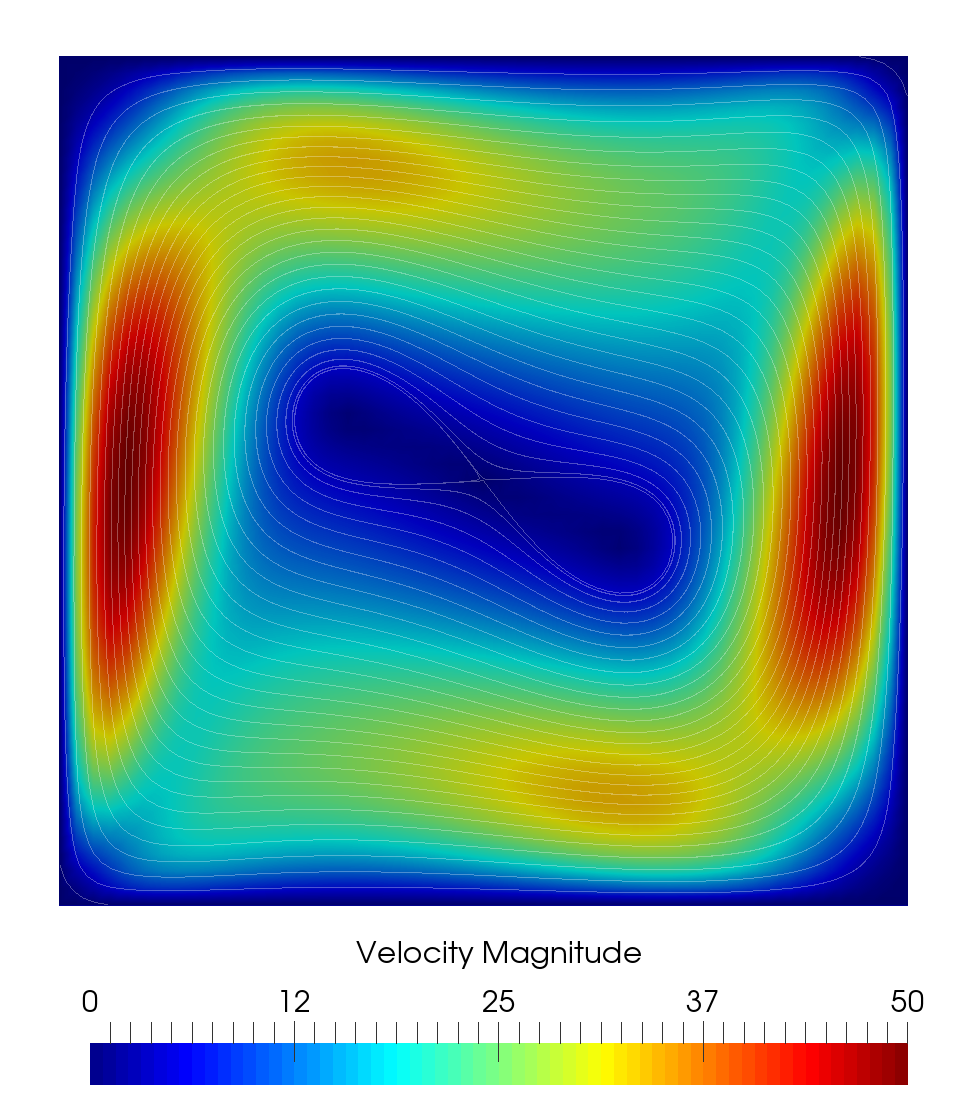}
\includegraphics[width=0.3\linewidth]{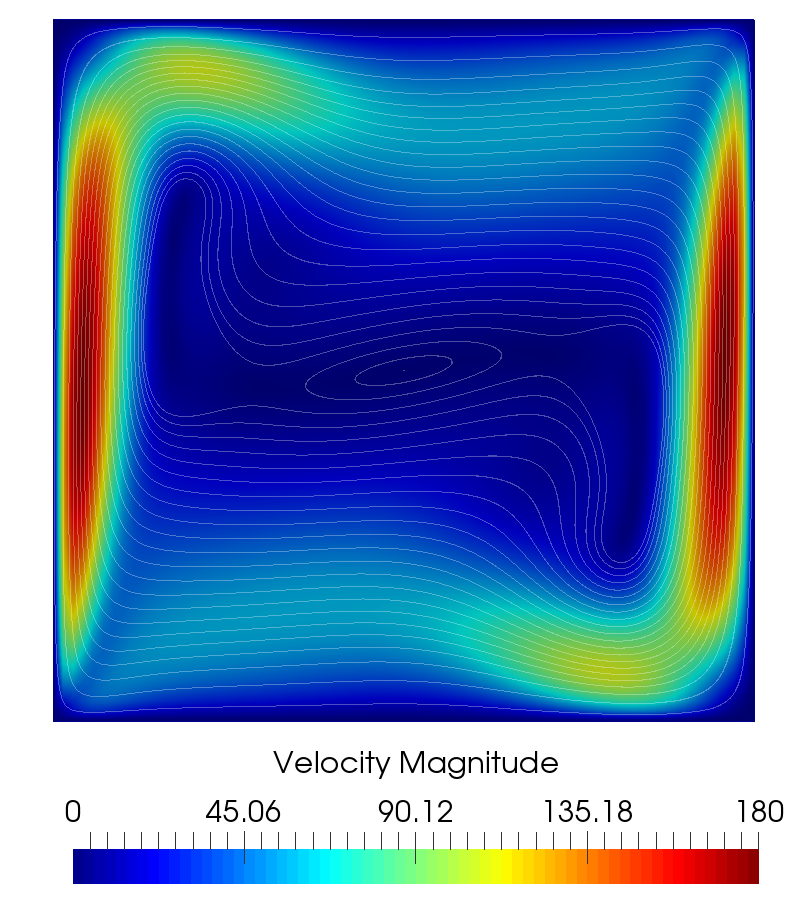}

\includegraphics[width=0.3\linewidth]{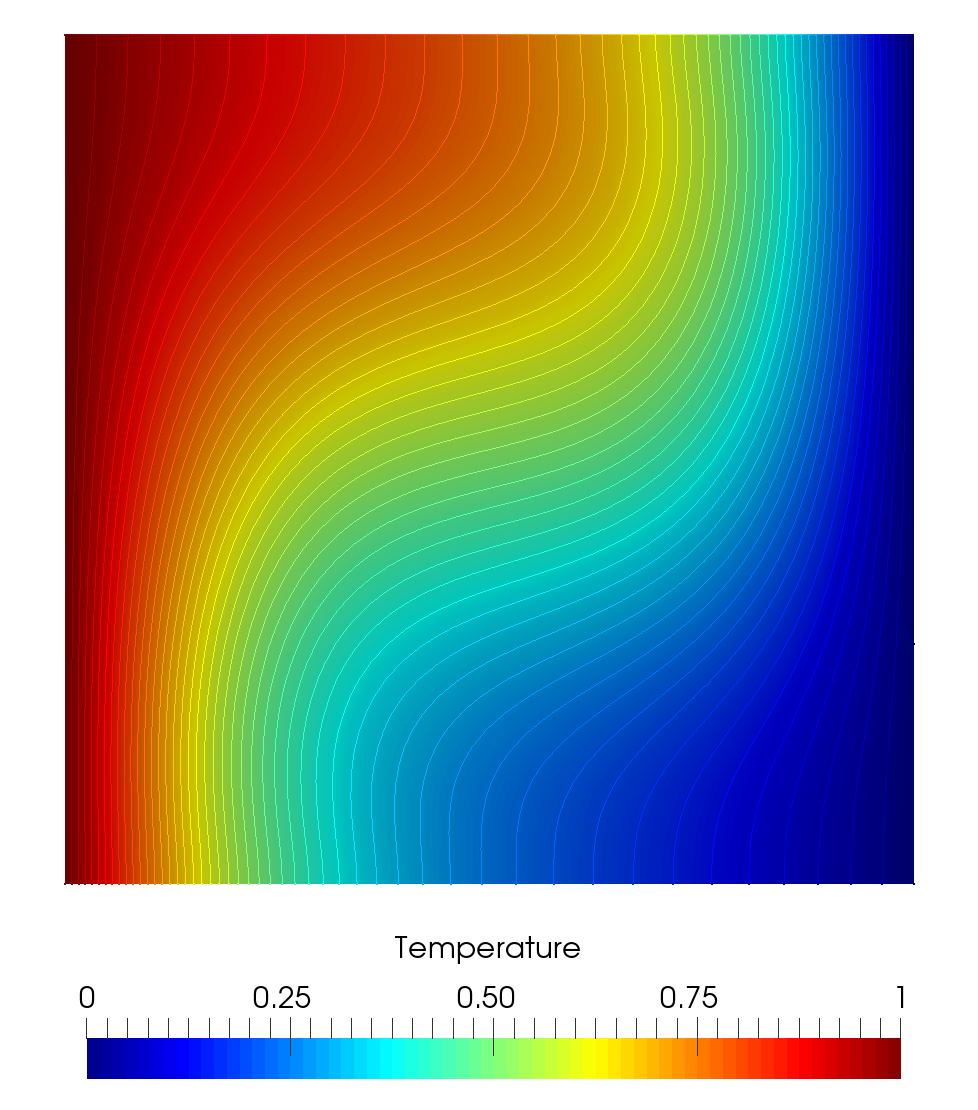}
\includegraphics[width=0.3\linewidth]{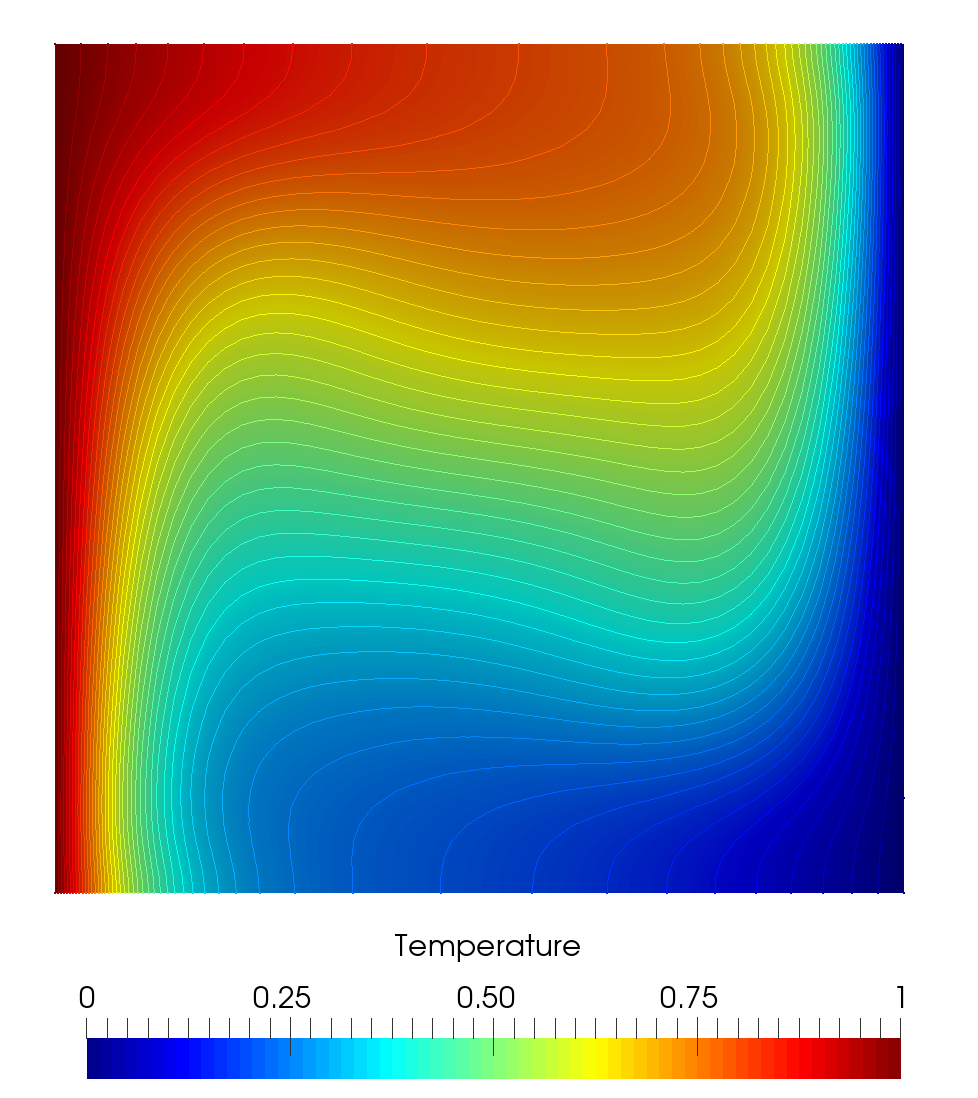}
\includegraphics[width=0.3\linewidth]{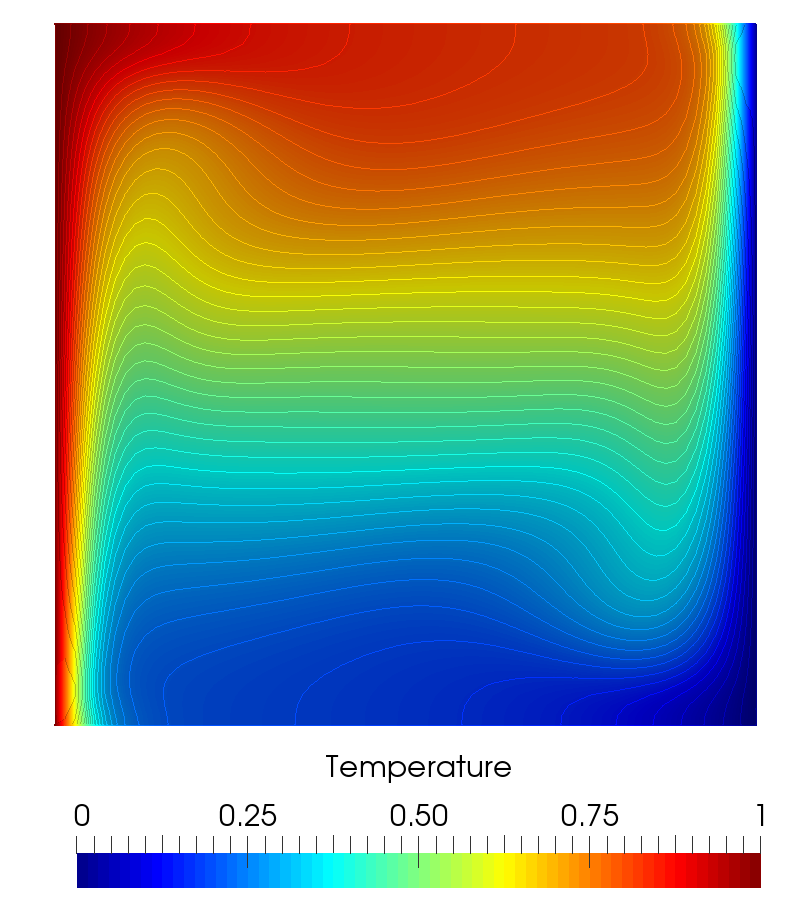}
\caption{FE solution, velocity magnitude (top) and temperature (bottom), for $\muf=4363$, $\muf=53778$ and  $\muf=667746$ (left to right), $\mug=1$. }\label{chap:Bouss::fig:Ra4363}
\end{figure}

In this test, we consider \black{two different scenarios}, one for the Rayleigh number range $[10^3, 10^5]$, an the other for the Rayleigh number range $[10^5, 10^6]$. In both cases we consider the geometrical parameter fixed, with $\mug=1$. For the first scenario, which corresponds to the lower Rayleigh number values, %(see e.g. \cite{Boussinesq_FE, square_cavity, Wan2001}) 
the heat transfer is principally in form of diffusion, i.e., the diffusion term in the energy equation is predominant, leading to an almost vertical linear contouring for the temperature, and a recirculating motion in the core of the region is observed. As we increase the value of the Rayleigh number in $\cD$, the flow is stretched to the walls, especially to the vertical walls; and the heat transfer \red{starts} to be driven mainly by convection. The isotherms become horizontal in a domain inside the cavity, far from the walls, that increases as the Rayleigh number increases. When we consider the second scenario, where the Rayleigh number range is higher, the velocity in the center of the cavity is practically zero, and presents large and normal gradients near the vertical walls. The temperature isolines are horizontal in a large domain inside the cavity, \black{except near the vertical walls}. This behavior agrees with the results presented in several works, e.g. \cite{Boussinesq_FE,square_cavity,Wan2001}. %This behavior can be observed in Figures \ref{chap:Bouss::fig:Ra4363} and \ref{chap:Bouss::fig:Ra53778}. 
In Fig. \ref{chap:Bouss::fig:Ra4363} we show the FE velocity magnitude and temperature for $\muf=4363$, $\muf=53778$ and $\muf=667746$, with a fixed value of the geometrical parameter of $\mug=1$.

\black{We consider different meshes depending the scenario}. For $\muf\in[10^3,10^5]$ we consider a uniform mesh, with 50 divisions in each square side, i.e., $h=0.02\sqrt{2}$. For $\muf\in[10^5,10^6]$ we consider a finer mesh, with 70 divisions in each square side, i.e., $h=1/70\cdot\sqrt{2}$, in order to reproduce efficiently the eddies near the vertical walls \black{appearing in this Rayleigh number range}.

Concerning the time step in the evolution semi-implicit approach, we have considered a time step $\Delta t=0.01$ for the case of $\muf\in[10^3,10^5]$, and $\Delta t=2\cdot10^{-3}$ for the case $\muf\in[10^5,10^6]$.

\begin{figure}[h!]
\centering
\includegraphics[width=0.49\linewidth]{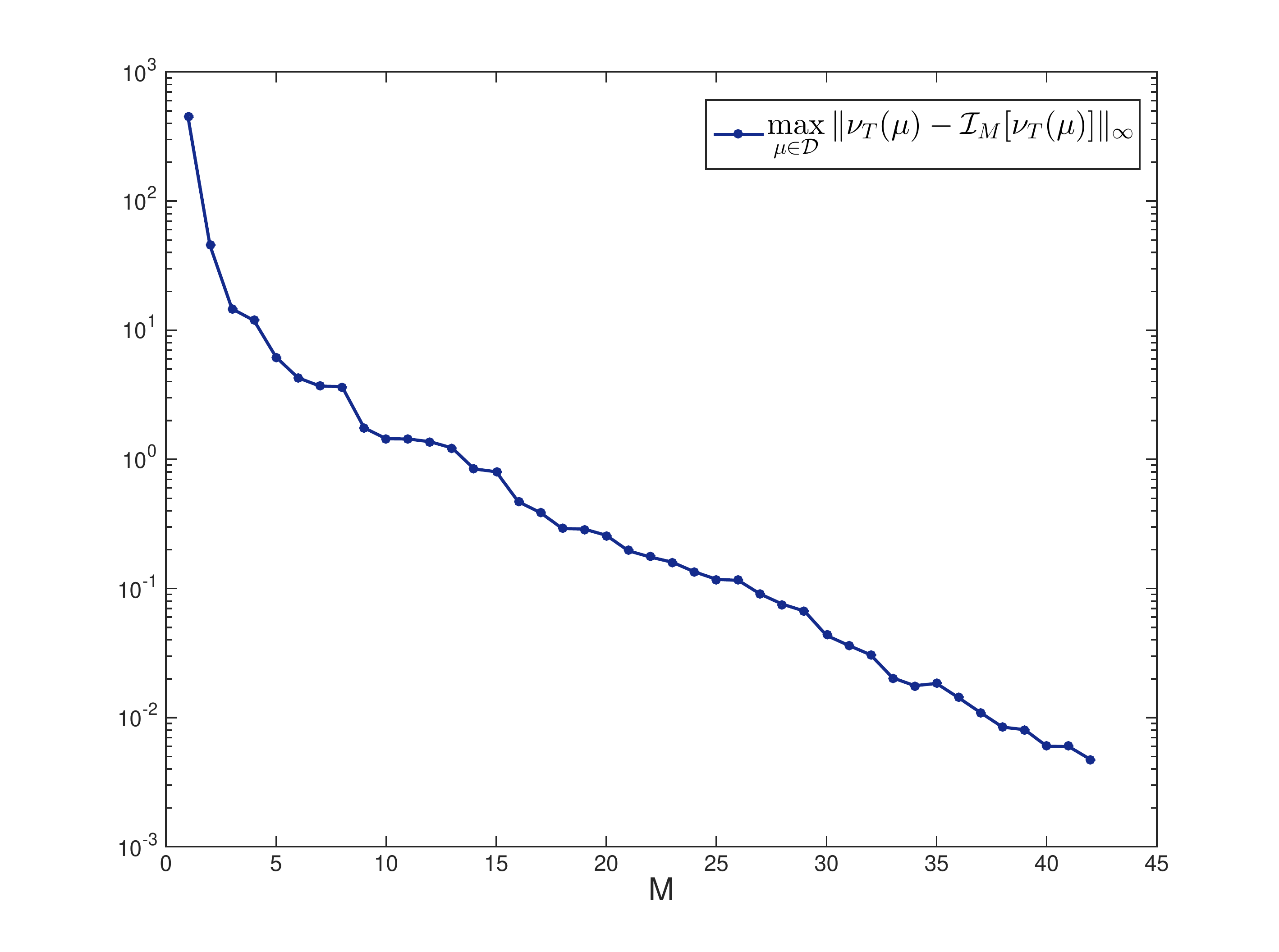}
\includegraphics[width=0.49\linewidth]{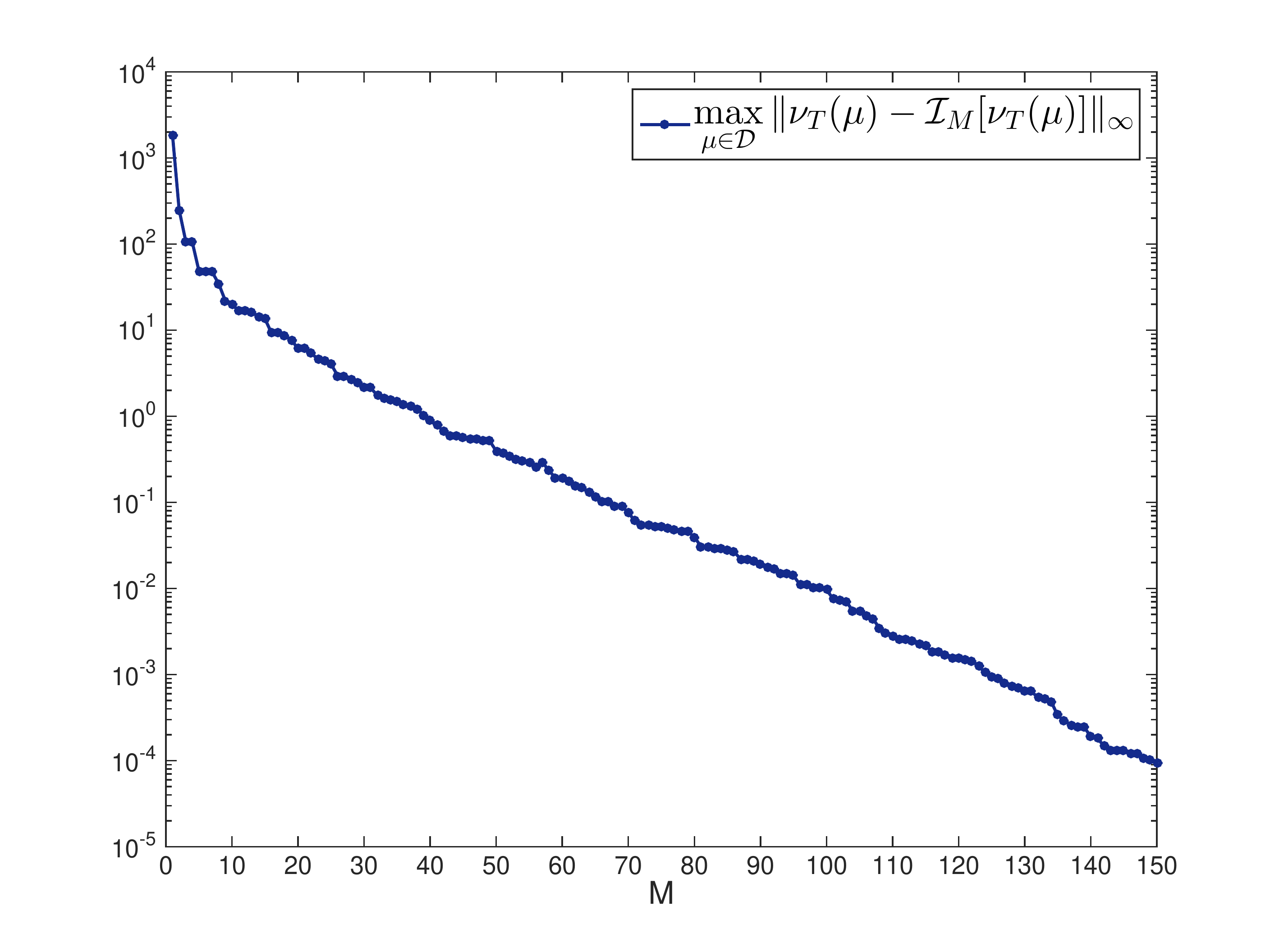}
\caption{Error evolution for the EIM, for $\muf\in[10^3,10^5]$ (left) and $\muf\in[10^5,10^6]$ (right).}\label{chap:Bouss::fig:EIM35}
\end{figure}

In the Reduced Basis framework, we perform an EIM for both the eddy viscosity and eddy diffusivity. Although for the numerical analysis performed in this work we have considered a regularized eddy diffusivity, the numerical tests are done with the eddy diffusivity defined in (\ref{chap:Bouss::eq:eddy_conduct}). Since the eddy diffusivity is proportional to the eddy viscosity, we only need to perform one EIM. With the EIM we are able to decouple the parameter dependence of the non-linear eddy viscosity and eddy diffusivity terms. For this test, we need $M=42$ basis until reaching a prescribed tolerance of $\veps_{EIM}=5\cdot10^{-3}$, when we consider that $\muf\in[10^3,10^5]$, and $M=150$ basis functions when we consider the second scenario where $\muf\in[10^5,10^6]$. In this last case, the Smagorinsky eddy viscosity and eddy diffusivity terms become more relevant, and for this reason, we take a lower tolerance for this test \black{with respect to the previous one}, considering $\veps_{EIM}=10^{-4}$. In Fig. \ref{chap:Bouss::fig:EIM35} we show the evolution of this error for both scenarios.

\begin{figure}[H]
\centering
\includegraphics[width=0.495\linewidth]{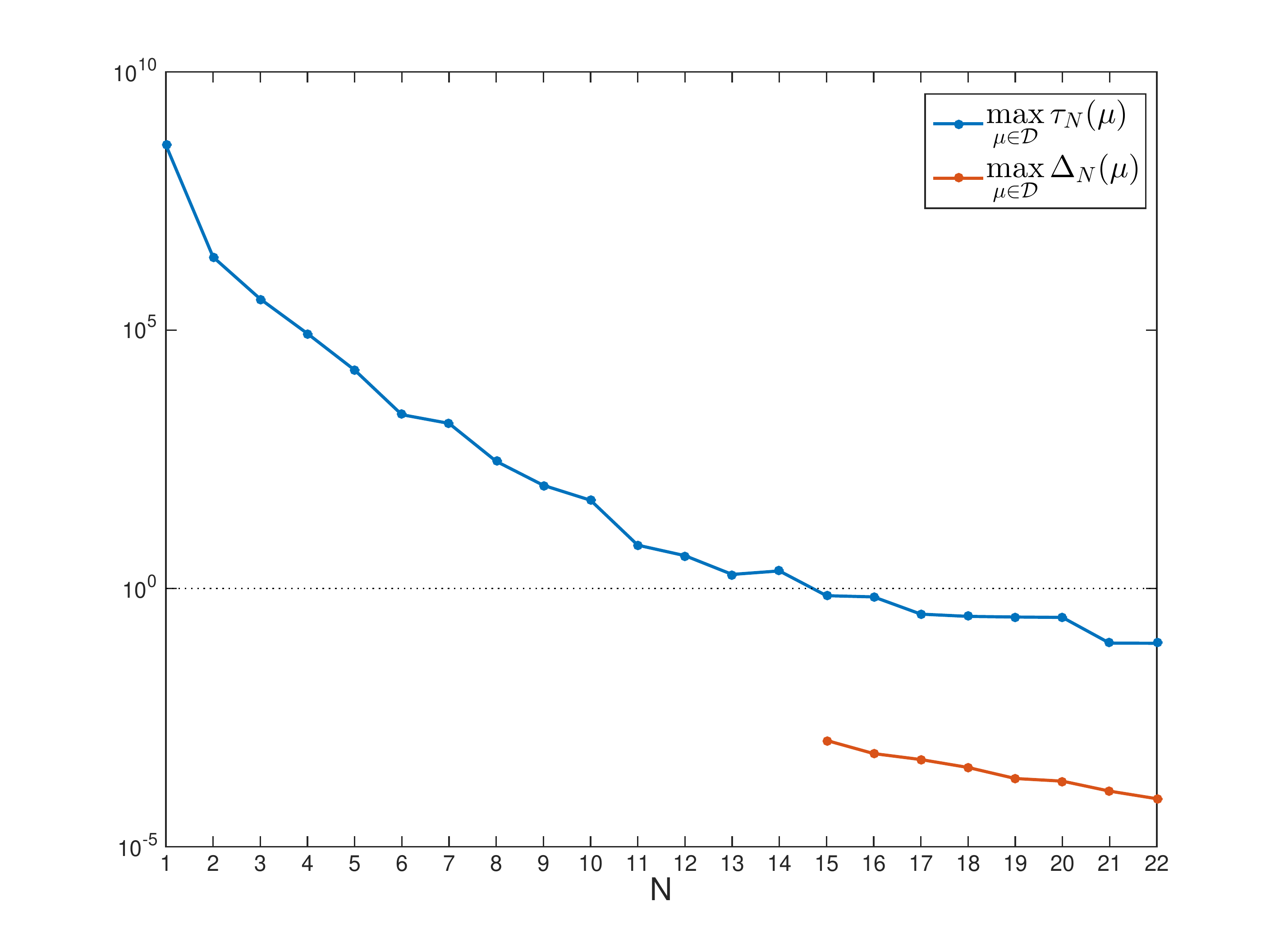}
\includegraphics[width=0.495\linewidth]{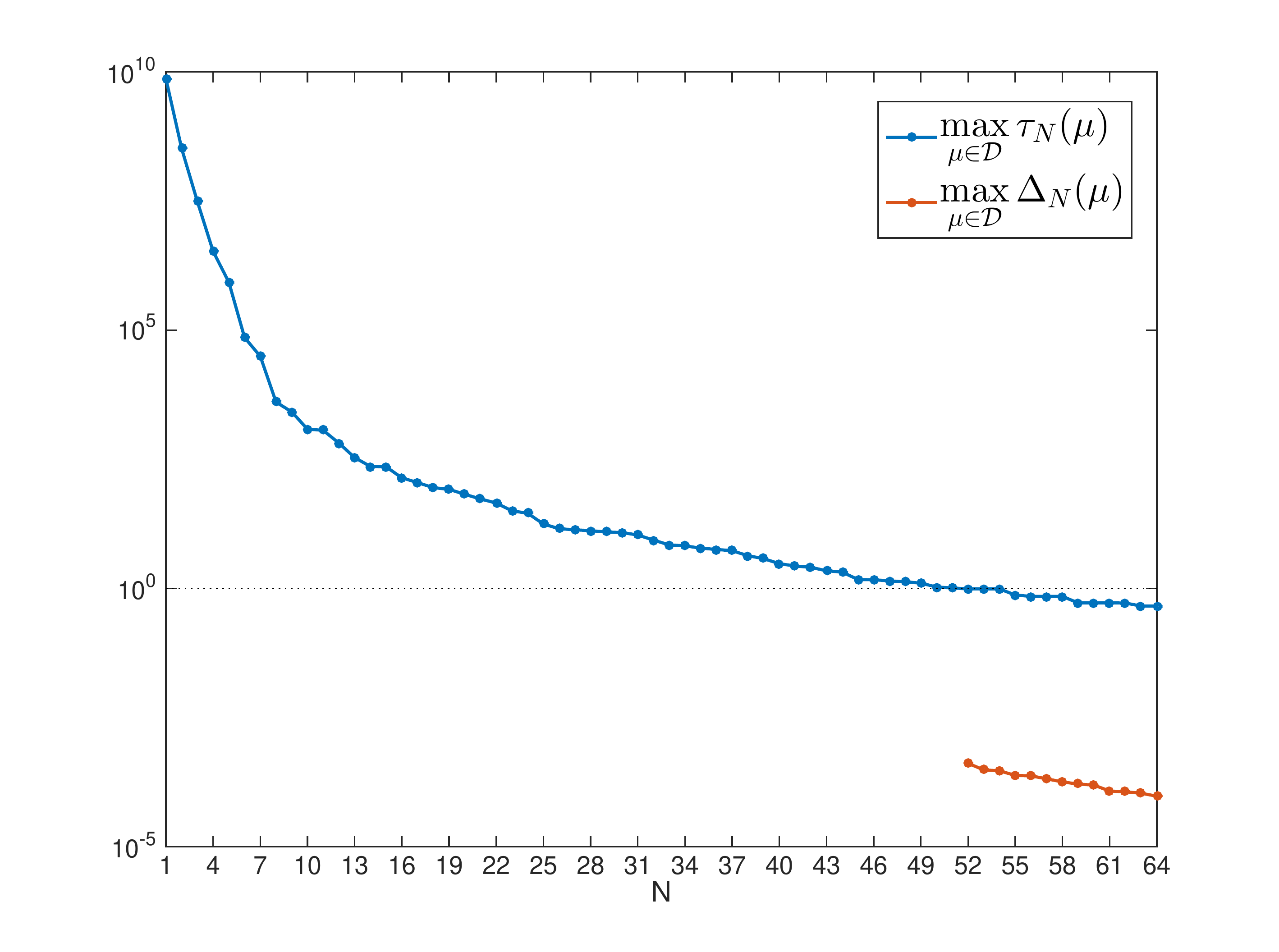}
\caption{Evolution of the \textit{a posteriori} error bound in the Greedy algorithm, for $\muf\in[10^3,10^5]$ (left) and $\muf\in[10^5,10^6]$ (right).}\label{chap:Bouss::fig:Delta35}
\end{figure}

For the Greedy algorithm we prescribe a tolerance of $\veps_{RB}=10^{-4}$ for both scenarios. For the first scenario, when $\muf\in[10^3,10^5]$, we need $N_{\max}=22$ basis to reach this tolerance. When $N=15$, holds the condition of Theorem \ref{chap:Geom::teor:Teorprinc} and $\tau_N(\muk)<1$ for all $\muk$ in $\cD$. In the second scenario, when $\muf\in[10^5,10^6]$, we need $N=N_{\max}=64$ basis functions to reach the tolerance previously prescribed, becoming $\taunk$ smaller than one when we get $N=52$ basis functions.
In both cases, when $\tau_N(\muk)>1$ and the \textit{a posteriori} error bound is not defined, we use as \textit{a posteriori} error bound the proper $\taunk$. In Fig. \ref{chap:Bouss::fig:Delta35} we show the convergence for the greedy algorithm.

\begin{figure}[H]
\centering
\includegraphics[width=0.495\linewidth]{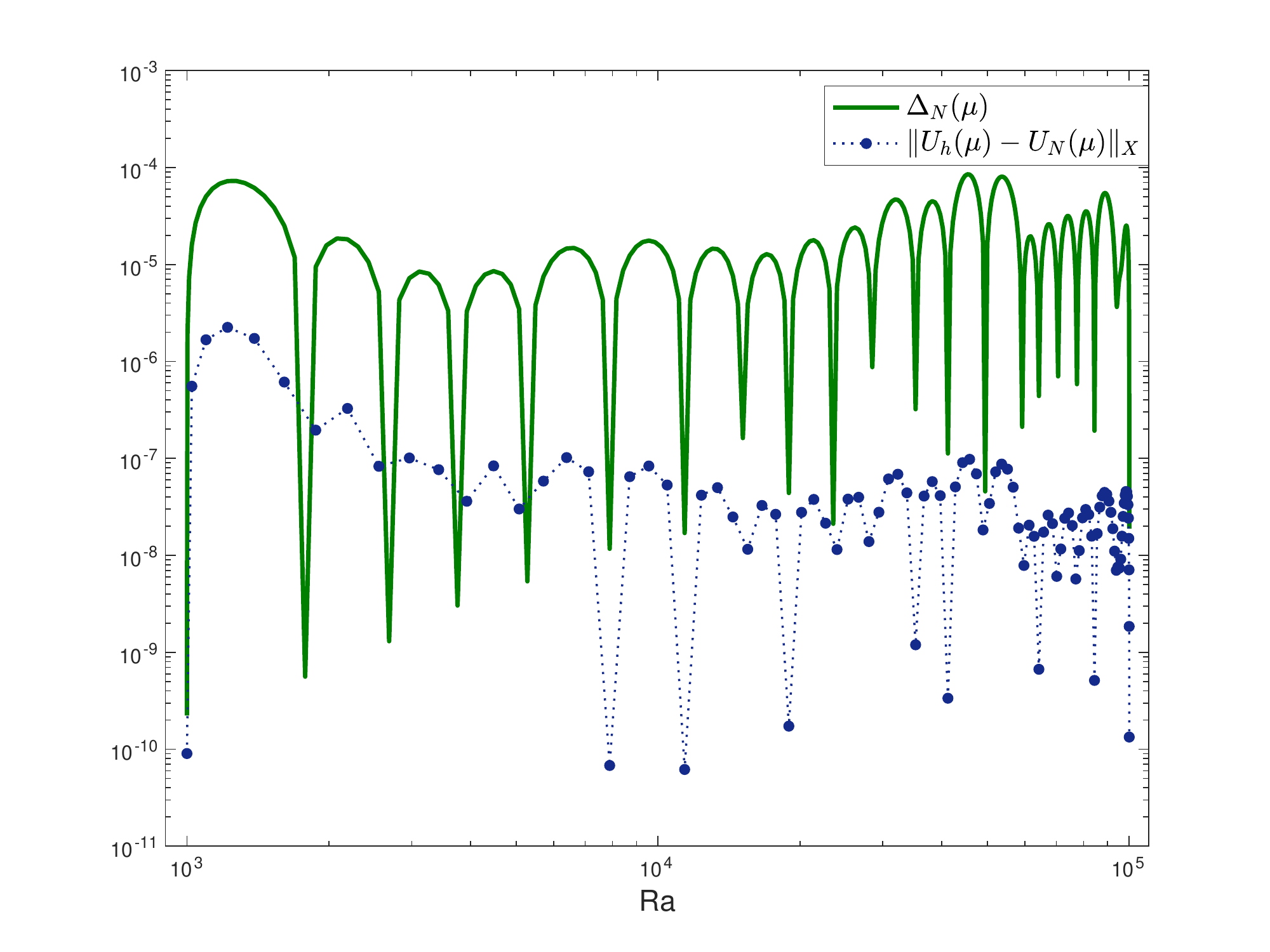}
\includegraphics[width=0.495\linewidth]{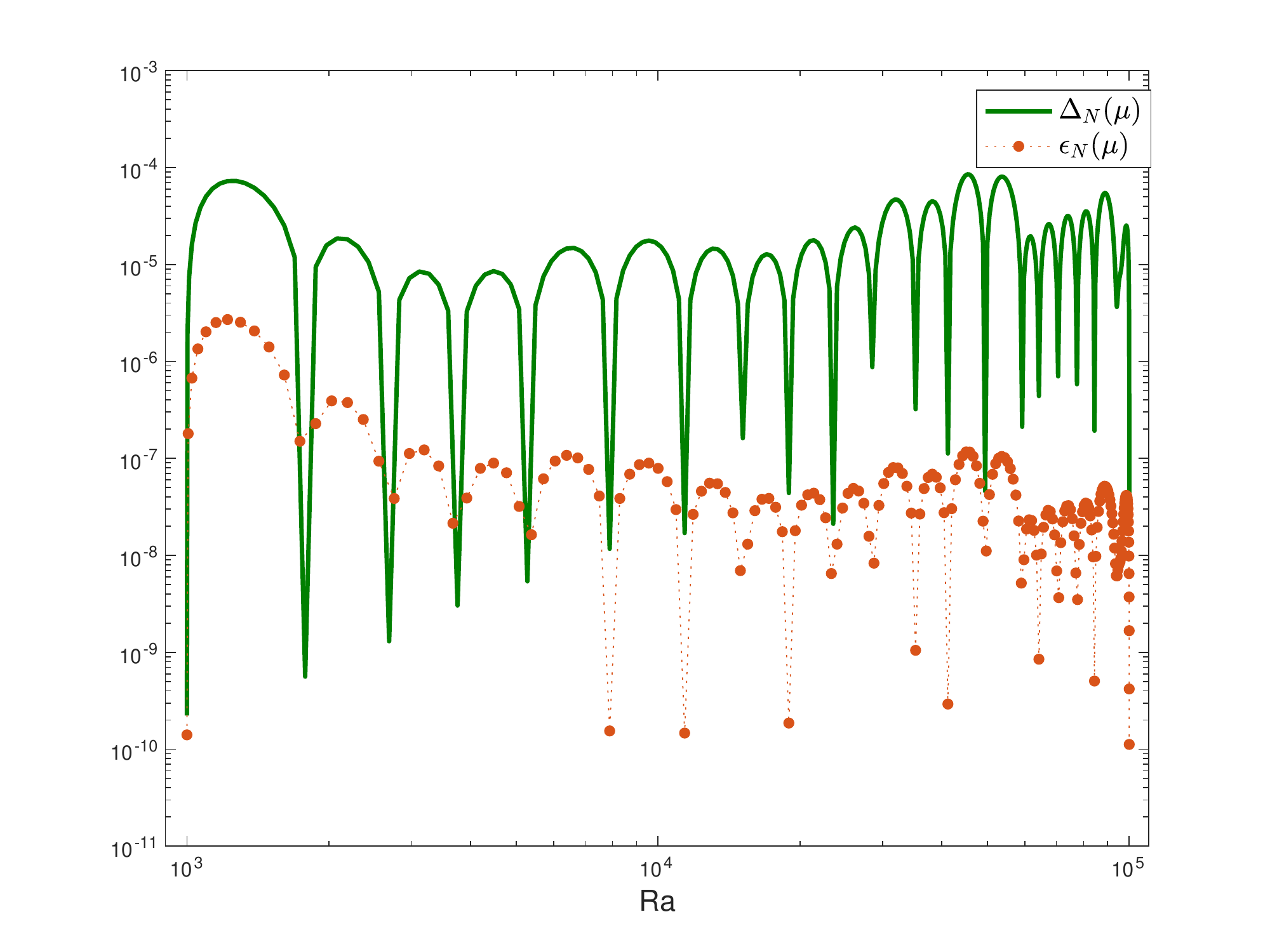}
\caption{\textit{A posteriori} error bound vs true error for $N=N_{\max}$ (left) and \textit{A posteriori} error bound vs $\en$ (right), for $\muf\in[10^3,10^5]$.}\label{chap:Bouss::fig:DeltaN35}
\end{figure}

\red{In Fig. \ref{chap:Bouss::fig:DeltaN35} (left) we show the comparison between the true error and the \textit{a posteriori} error bound when $\muf\in[10^3,10^5]$, in which we can observe that the efficiency of the \textit{a posteriori} error bound is between two and three orders of magnitude. Moreover, in Fig. \ref{chap:Bouss::fig:DeltaN35} (right) we represent the comparison between the true error and the dual norm of the residual, $\en$. For this test, we can observe how the error correlates quite better with the dual norm of the residual than the \textit{a posteriori} error bound, but in any case, it does not give us an upper bound of the error.}

\begin{figure}[H]
\centering
\includegraphics[width=0.495\linewidth]{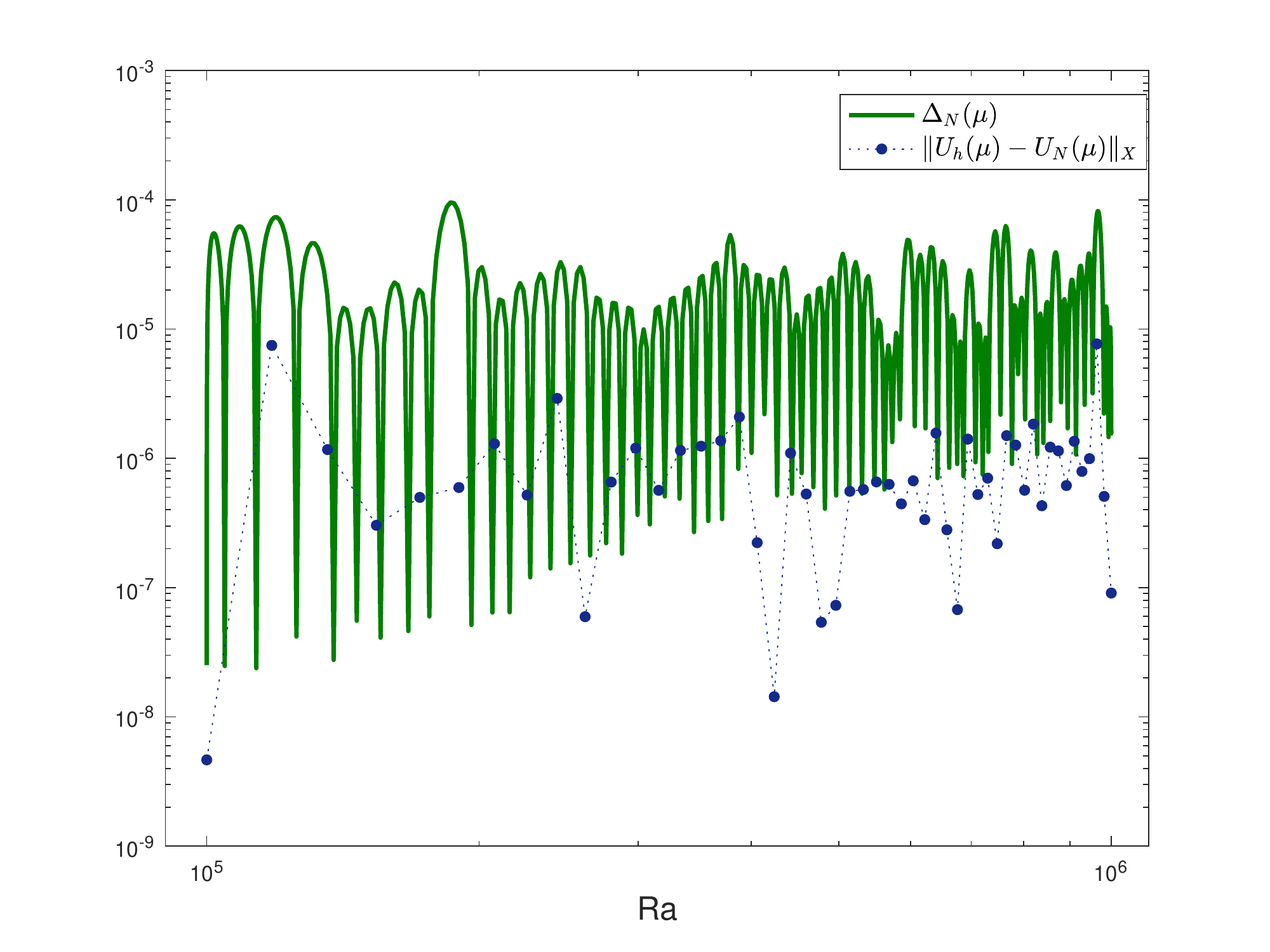}
\includegraphics[width=0.495\linewidth]{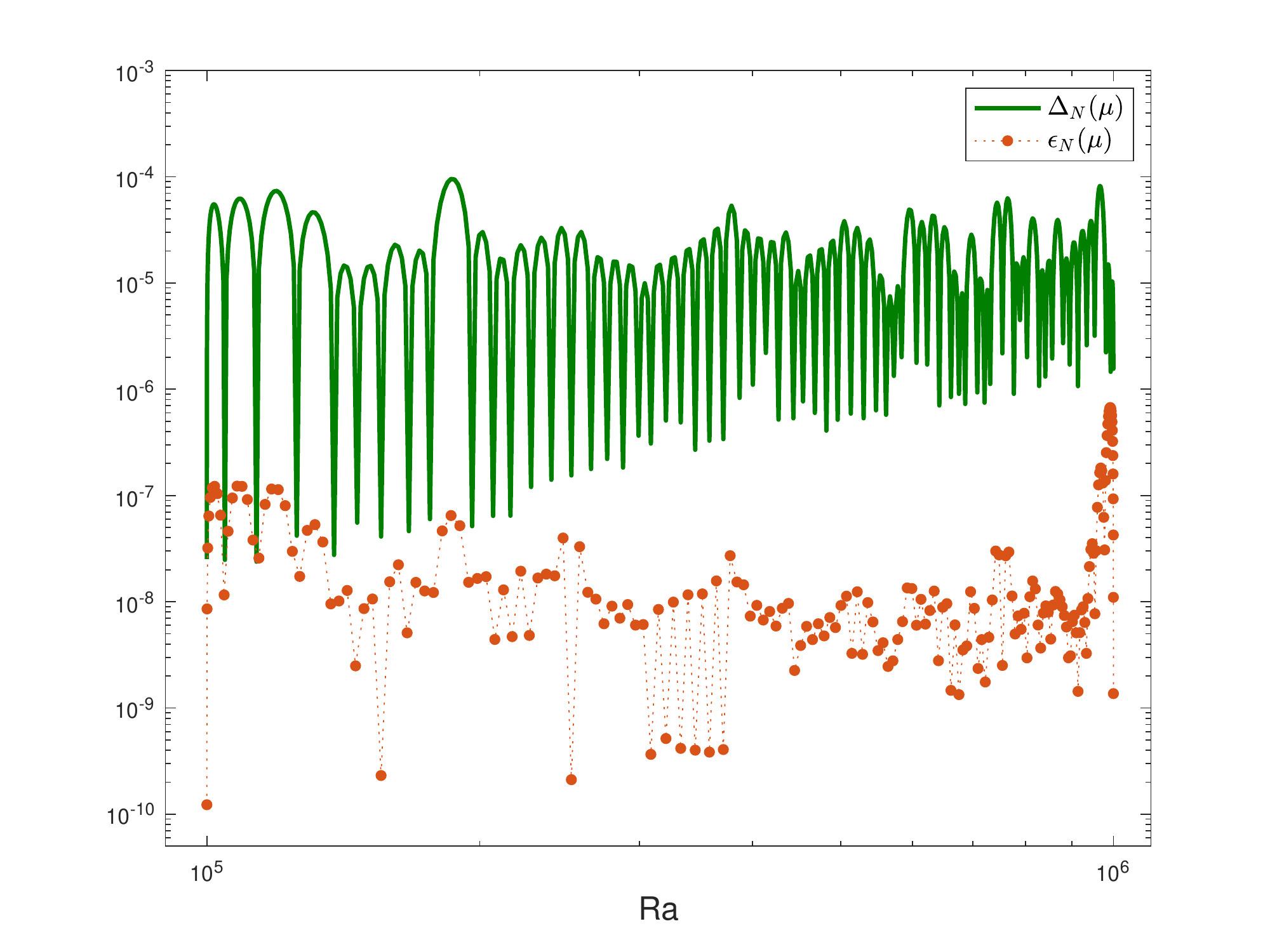}
\caption{\textit{A posteriori} error bound vs true error for $N=N_{\max}$ (left) and \textit{A posteriori} error bound vs $\en$ (right), for $\muf\in[10^5,10^6]$.}\label{chap:Bouss::fig:DeltaN56}
\end{figure}

\red{In Fig. \ref{chap:Bouss::fig:DeltaN56} (left) we show the comparison between the true error and the \textit{a posteriori} error bound when $\muf\in[10^5,10^6]$, in which we can observe that the efficiency of the \textit{a posteriori} error bound is between one and two orders of magnitude. Moreover, in Fig. \ref{chap:Bouss::fig:DeltaN56} (right) we represent the comparison between the true error and the dual norm of the residual, $\en$. For this test case, we can see how the dual norm of the residual is much lower than the error. Thus, the \textit{a posteriori} error estimator gives us a better estimation for the error than the dual norm of the residual.}

\medskip
\medskip

Finally, in Table \ref{chap:Bouss::tab:Ra_bajo}, we show a comparison between the FE and RB solutions for several Rayleigh values in both scenarios. We show the computational time for solving a FE solution and a RB solution in the online phase. As can be observed, the speed-up rate of the computational time is \black{larger than three orders of magnitude} when $\mug\in[10^3,10^5]$, while when $\mug\in[10^5,10^6]$ the speed-up rate is \black{close to three hundred}. \black{The difference in the speed-up magnitude between both cases is due to the longer number of EIM and RB functions computed in each case}. In addition, we show the \red{relative} errors in $H^1$-norm for velocity and temperature, and in $L^2$-norm for pressure; for which we observe that the RB solution is close enough to the FE solution, \red{with relative errors  \red{about $10^{-9}$} in both test cases}. \black{For this test, the offline phase when $\muf\in[10^3,10^5]$ took approximately 2 days in being performed. For the case when $\muf\in[10^5,10^6]$, the offline phase took approximately 3 weeks in be performed. In this offline computational time we consider either the EIM and the Greedy algorithm with the computation of the a posteriori error estimator}

\begin{table}[h!]
$$
\begin{tabular}{l|cccc}
\hline
Data &$\muf=4060\quad$&$\muf=17808\;\;$&$\muf=53778\;\;$&$\muf=93692$\\
\hline
$T_{FE}$&633.65s & 585.83s & 553.25s&677.86s\\
$T_{online}$& 0.55s& 0.5s& 0.46s&0.49s\\
%\hline
%Iter FE& 34 & 72  & 135  &204  \\
%Iter RB& 25 & 52  & 88 & 112 \\
\hline
speedup& 1133 & 1151& 1189 &1367\\
\hline
$\dfrac{\|\uk_h-\uk_N\|_1}{\|\uk_h\|_1}$&$4.01\cdot10^{-9}$&$4.06\cdot10^{-9}$ & $3.56\cdot10^{-9}$ &$3.22\cdot10^{-9}$\\
\hline
$\dfrac{\|\theta_h-\theta_N\|_1}{\|\theta_h\|_1}$&$5.59\cdot10^{-9}$&$4.29\cdot10^{-9}$ & $5.36\cdot10^{-9}$ &$4.76\cdot10^{-9}$\\
\hline
$\dfrac{\|p_h-p_N\|_0}{\|p_h\|_0}$&$5.75\cdot10^{-10}$&$1.71\cdot10^{-10}$ & $1.96\cdot10^{-10}$  &$2.30\cdot10^{-10}$\\
\hline
\end{tabular}$$
%\hspace{-2cm}\caption{Computational time for FE and RB solutions, with the speedup and the error, for problem (\ref{chap:Geom::pb:Cont_Bouss}), and $Ra\in[10^3,10^5]$, $\mug=1$.}
\label{chap:Bouss::tab:Ra_bajo}
$$\begin{tabular}{l|cccc}
\hline
Data &$\muf=169411$&$\muf=355402$&$\muf=667746$&$\muf=921441$\\
\hline
$T_{FE}$&3563.11s & 3675.01s & 4354.26s&4928.37s\\
$T_{online}$& 9.28s& 11.34s& 15.22s&16.8s\\
\hline
speedup& 383 & 324& 285 & 293\\
\hline
$\dfrac{\|\uk_h-\uk_N\|_1}{\|\uk_h\|_1}$&$2.88\cdot10^{-9}$&$2.87\cdot10^{-9}$ & $1.59\cdot10^{-9}$ &$3.05\cdot10^{-9}$\\
\hline
$\dfrac{\|\theta_h-\theta_N\|_1}{\|\theta_h\|_1}$&$5.93\cdot10^{-9}$&$5.63\cdot10^{-9}$ & $5.11\cdot10^{-9}$ &$4.81\cdot10^{-9}$\\
\hline
$\dfrac{\|p_h-p_N\|_0}{\|p_h\|_0}$&$3.53\cdot10^{-10}$&$3.83\cdot10^{-10}$ & $3.79\cdot10^{-10}$  &$3.65\cdot10^{-10}$\\
\hline
\end{tabular}$$
\caption{Computational time for FE and RB solutions, with the speedup and the error, for \black{problem (\ref{chap:Geom::pb:Cont_Bouss}), $Ra\in[10^3,10^5]$} (top) and $Ra\in[10^5,10^6]$ (bottom), $\mug=1$.}\label{chap:Bouss::tab:Ra_alto}
\end{table}

\subsection{Geometrical parametrization}\label{chap:Geom::sec:geom_param}
\begin{figure}[h]
\centering
\includegraphics[width=0.8\linewidth]{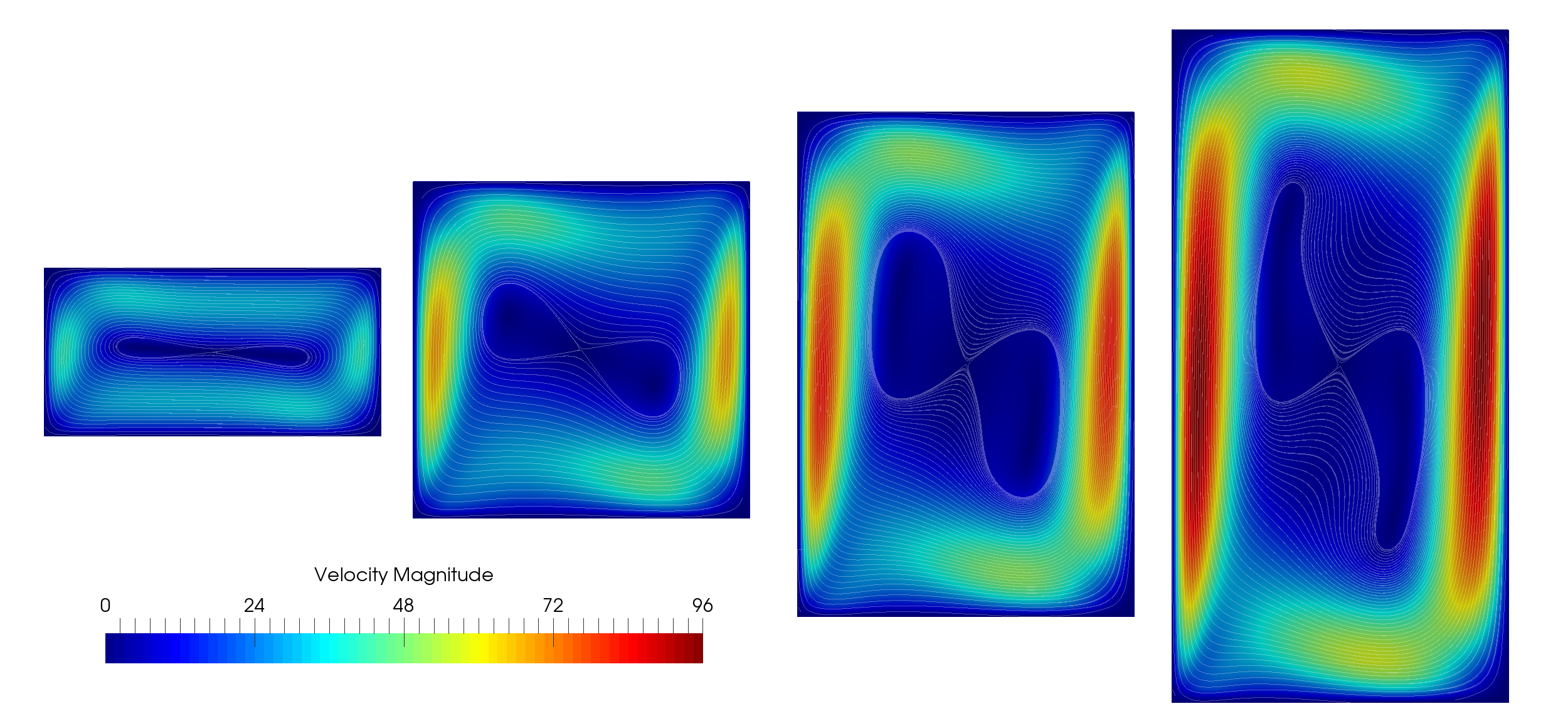}
\caption{FE snapshots for $\mug=0.5$, $\mug=1$, $\mug=1.5$ and $\mug=2$ (left to right).}\label{chap:Geom::fig:Geom_FE}
\end{figure}

In this test, we consider a moderate Rayleigh number value $Ra=10^5$, and we consider the geometrical parameter ranging in $\mug\in\cD=[0.5,2]$. The difference in the height of the cavity affects to the buoyancy force, making it more relevant when we increase the parameter value. This behavior is observed in Fig. \ref{chap:Geom::fig:Geom_FE}, in which we show four solutions for different values of the geometrical parameter.

\begin{figure}[h]
\centering
\includegraphics[width=0.49\linewidth]{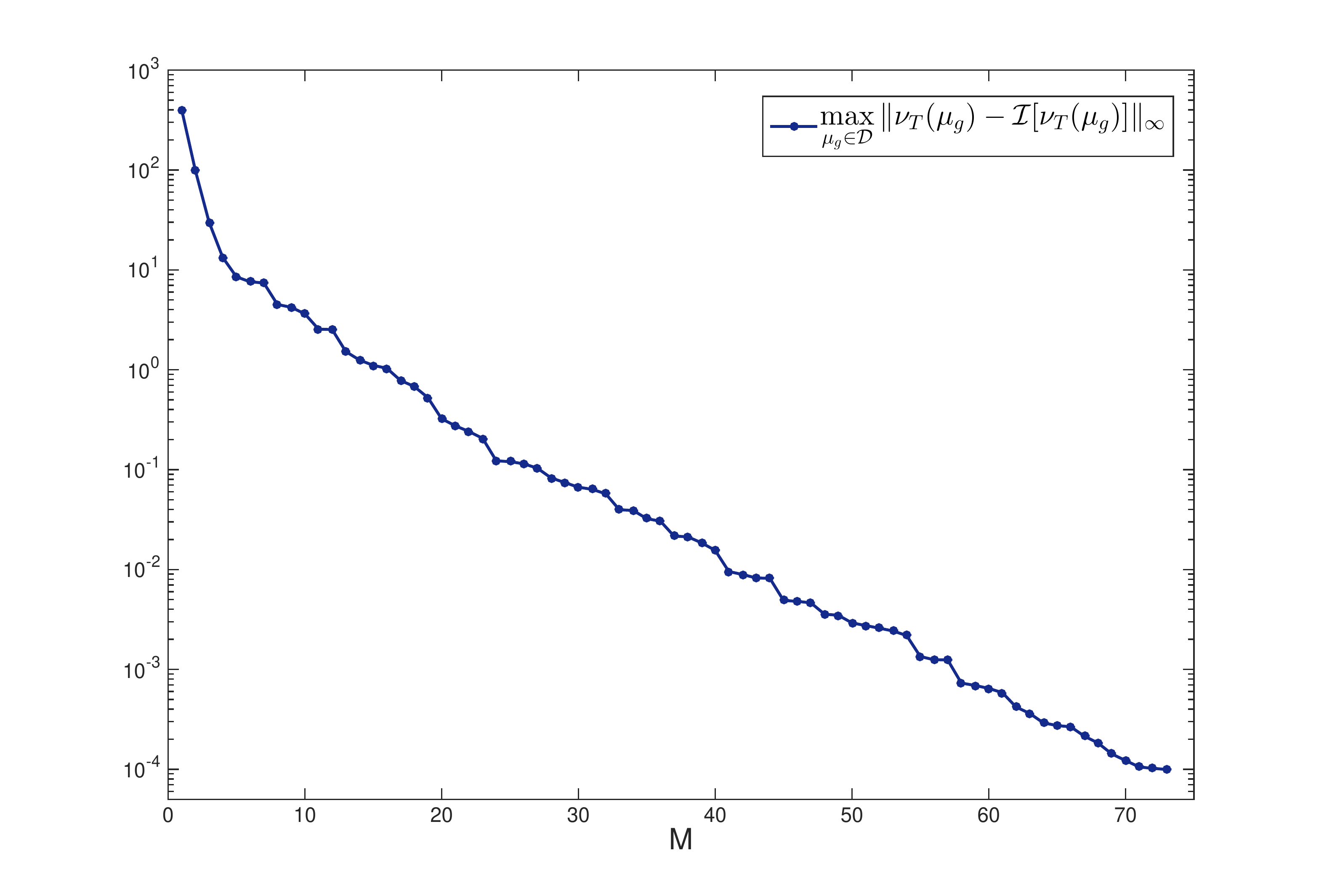}
\includegraphics[width=0.49\linewidth]{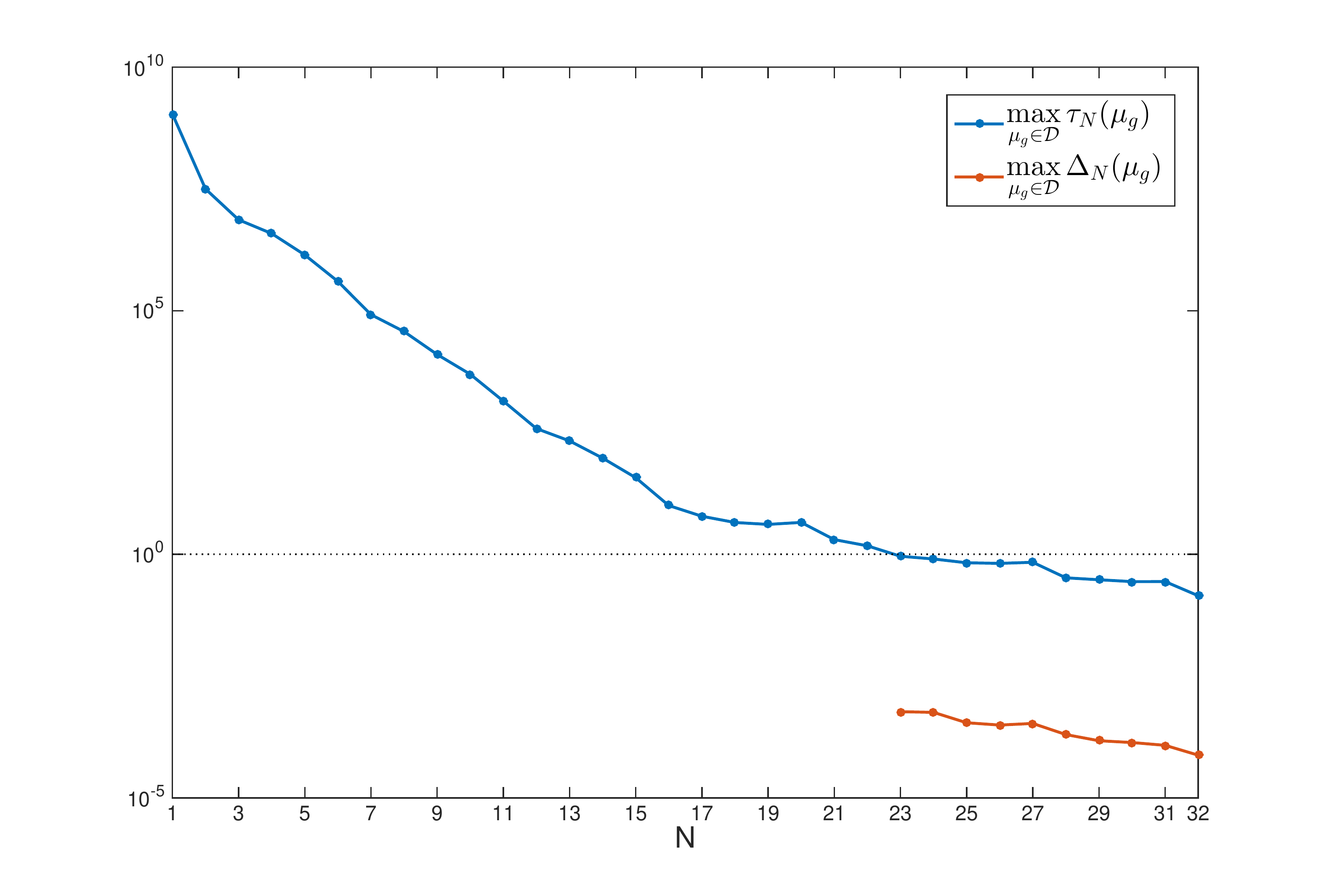}
\caption{Error evolution for the EIM (left) and Evolution of the \textit{a posteriori} error bound in the Greedy algorithm (right), for Boussinesq VMS-Smagorinsky model with $\mug\in[0.5,2]$.}\label{chap:Geom::fig:EIM}
\end{figure}

Firstly in the offline phase, we construct the reduced-basis space corresponding to the EIM, in which we approximate properly the \black{eddy viscosity and eddy diffusivity terms}. In this test, we need $M=73$ basis functions in order to reach a \black{prescribed tolerance} of $\veps_{EIM}=10^{-4}$. In Fig. \ref{chap:Geom::fig:EIM} (left) we show the evolution of the infinity norm of the error between the eddy viscosity $\nu_T(\mu_g)$ and its EIM approximation.

%\begin{figure}[h]
%\centering
%
%\caption{Evolution of the \textit{a posteriori} error bound in the Greedy algorithm, for Boussinesq VMS-Smagorinsky model with $\mug\in[0.5,2]$.}\label{chap:Geom::fig:Delta}
%\end{figure} 

For the Greedy algorithm, in this test, we prescribe a tolerance for the \textit{a posteriori} error bound of $\veps_{RB}=10^{-4}$. We need $N=23$ basis functions until \black{to guarantee the condition} of Theorem \ref{chap:Geom::teor:Teorprinc}, and get $\taunk<1$. Then, we reach the prescribed tolerance when $N=N_{\max}=32$. In Fig. \ref{chap:Geom::fig:EIM} (right), we show the maximum value for all $\mug\in\cD$ of the \textit{a posteriori} error estimator, and $\taunk$, in each iteration of the Greedy algorithm. 

\begin{figure}[h]
\centering
\includegraphics[width=0.49\linewidth]{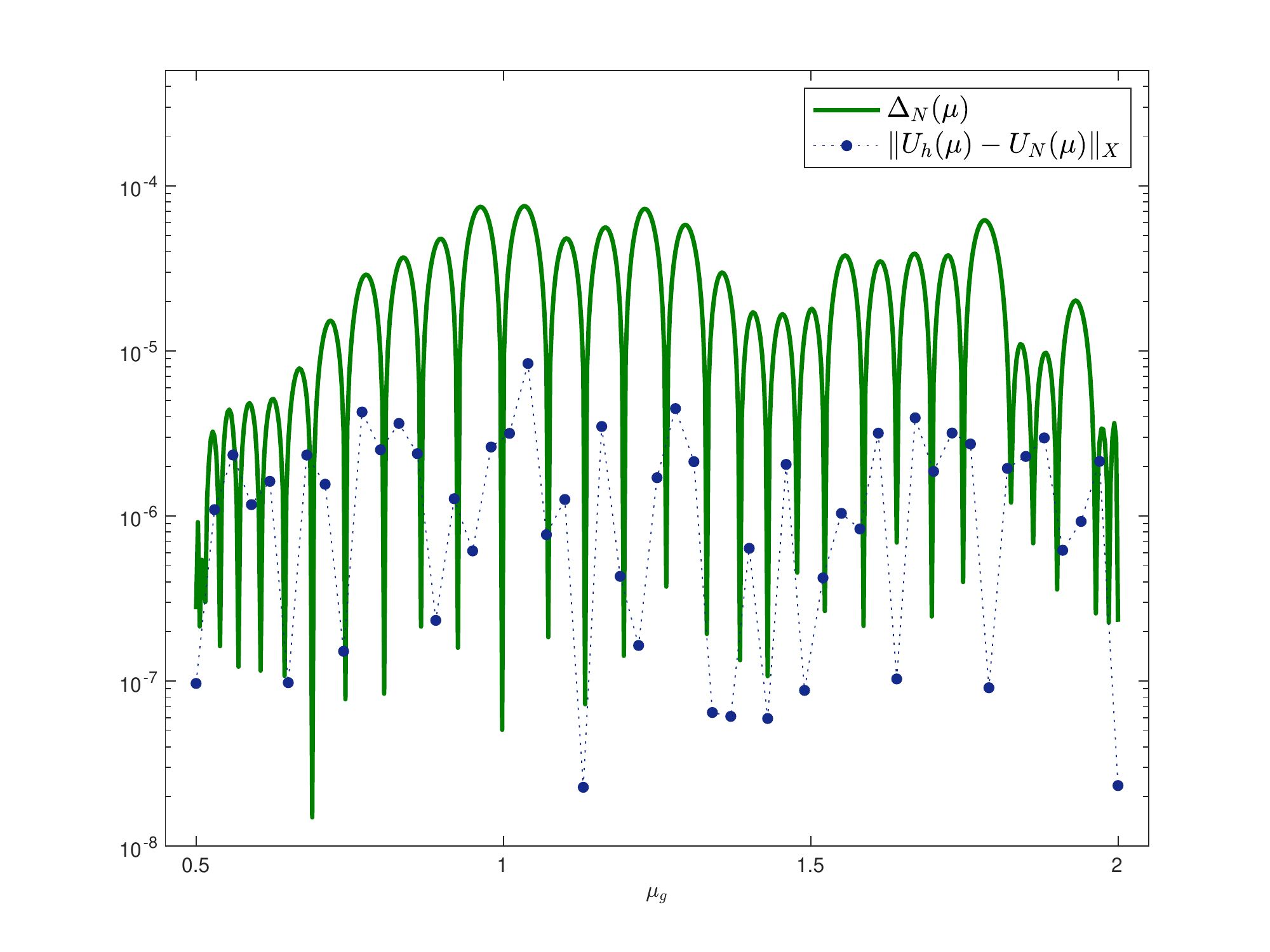}
\includegraphics[width=0.49\linewidth]{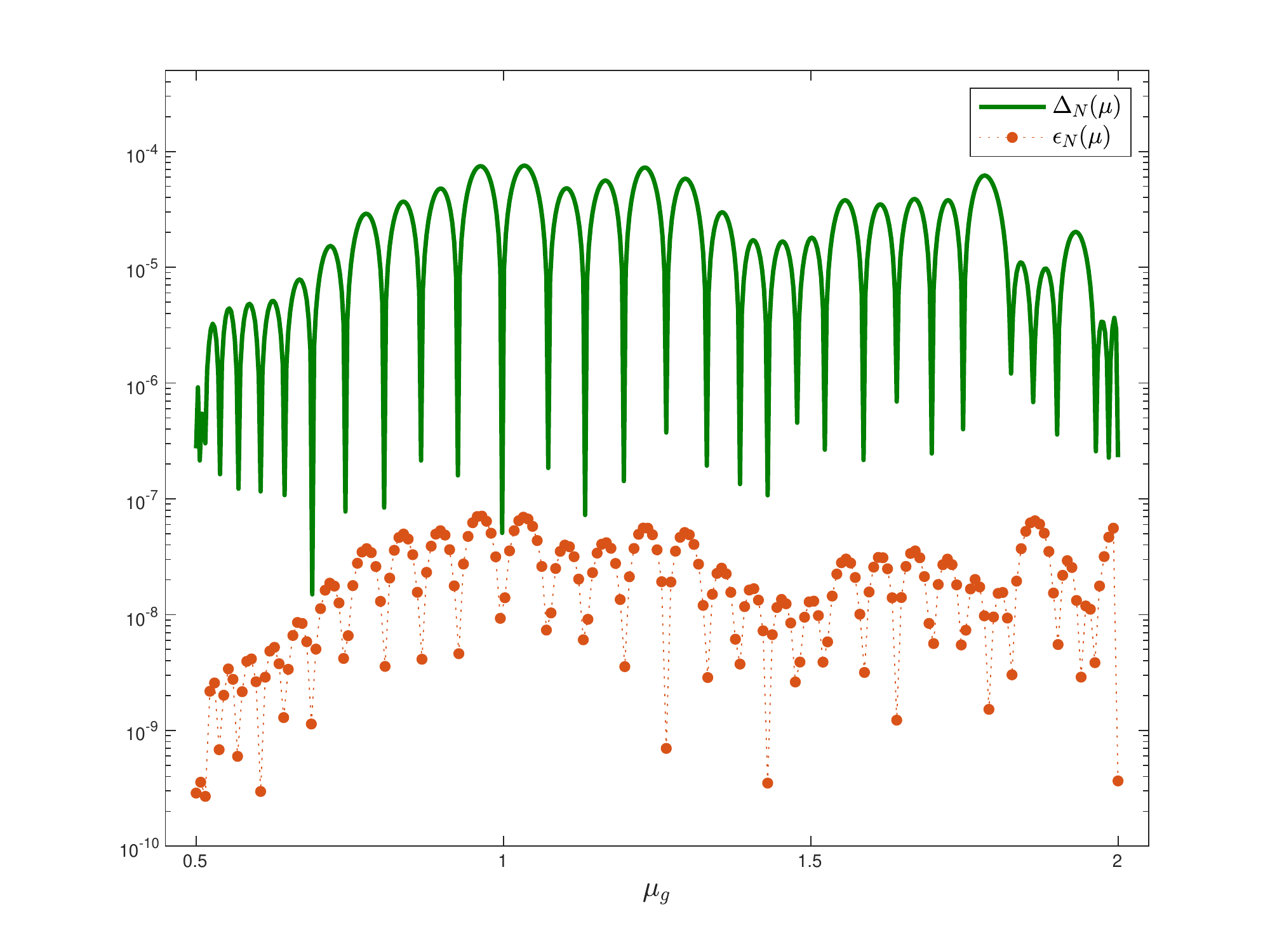}
\caption{\textit{A posteriori} error bound for $N=N_{\max}=32$.}\label{chap:Geom::fig:DeltaN}
\end{figure}

\red{In Fig. \ref{chap:Geom::fig:DeltaN} (left), we show a comparison between the \textit{a posteriori } error bound and the true error for all $\mug\in\cD$, in the last iteration of the Greedy algorithm, i.e., when $N=32$. Here we can see how the efficiency of the \textit{a posteriori} error bound is about one order of magnitude. In Fig. \ref{chap:Geom::fig:DeltaN} (right), we show the comparison between the \textit{a posteriori} error bound and $\en$. In this case, we can observe how the dual norm of the residual is much lower than the error, and how the \textit{a posteriori} error bound fits better the true error than the dual norm of the residual.}

Finally, in Table \ref{chap:Geom::tab:geom}, we summarize the results for several parameter values. We show the comparison between the time for computing a FE solution, and the online phase computational time. We obtain a speed-up rate of several hundreds in the computational time. The RB solution accuracy is fairly good, since the \red{relative error is approximately of order $10^{-9}$ for velocity, temperature, and pressure}. \black{For this test, the offline phase took approximately 5 days in being performed. In this offline computational time we consider either the EIM and the Greedy algorithm with the computation of the a posteriori error estimator}

\begin{table}[H]
$$
\begin{tabular}{l|cccc}
\hline
Data &$\mu_g=0.64$&$\mu_g=1.08$&$\mug=1.44$&$\mug=1.87$\\
\hline
$T_{FE}$&808.91s & 810.16s & 866.1s&851.82s\\
$T_{online}$& 2.68s& 2.55s& 2.61s&2.52s\\
%\hline
%Iter FE& 34 & 72  & 135  &204  \\
%Iter RB& 25 & 52  & 88 & 112 \\
\hline
speedup& 301 & 317& 331 &337\\
\hline
$\|\uk_h-\uk_N\|_1/\|\uk_h\|_1$&$3.4\cdot10^{-9}$&$4.12\cdot10^{-9}$ & $5.41\cdot10^{-9}$ &$5.68\cdot10^{-9}$\\
\hline
$\|\theta_h-\theta_N\|_1/\|\theta_h\|_1$&$3.75\cdot10^{-8}$&$4.66\cdot10^{-9}$ & $4.86\cdot10^{-9}$ &$4.91\cdot10^{-9}$\\
\hline
$\|p_h-p_N\|_0/\|p_h\|_0$&$2.51\cdot10^{-9}$&$3.25\cdot10^{-9}$ & $5.51\cdot10^{-9}$  &$4.48\cdot10^{-9}$\\
\hline
\end{tabular}$$
\hspace{-2cm}\caption{Computational time for FE and RB solutions, with the speed-up and the error, for Boussinesq VMS-Smagorinsky model with $\mug\in[0.5,2]$.}
\label{chap:Geom::tab:geom}
\end{table}

\subsection{Physical and geometrical parametrization}

In this test, we perform a RB model in which a physical parameter (the Rayleigh number), and a geometric parameter are taken into account. Due to the increasing complexity in the flux with the consideration of this two parameters, we consider low range of Rayleigh number. 

Thus, we consider that $\muk=(\muf,\mug)\in\cD=[10^3,10^4]\times[0.5,2]$. If we wanted to increase the Rayleigh number, we would have to consider a smaller interval for the geometric parameter. Indeed, as shown in sect. \ref{chap:Bouss::sec:High_Ra}, the flow for high Rayleigh values is quite complex, thus the consideration of geometric parameter joint with the physical parameter is only possible if both intervals are \black{not too big}. If a big parameter set is required, a possible strategy is to split it in subsets of \black{smaller} amplitude.%, as done in section \ref{chap:Bouss::sec:num}, where the whole interval of Rayleigh number of interest $[10^3,10^6]$ is split in two sub-intervals, $[10^3,10^5]$ and $[10^5,10^6]$. 

\begin{figure}[H]
\centering
\includegraphics[width=0.49\linewidth]{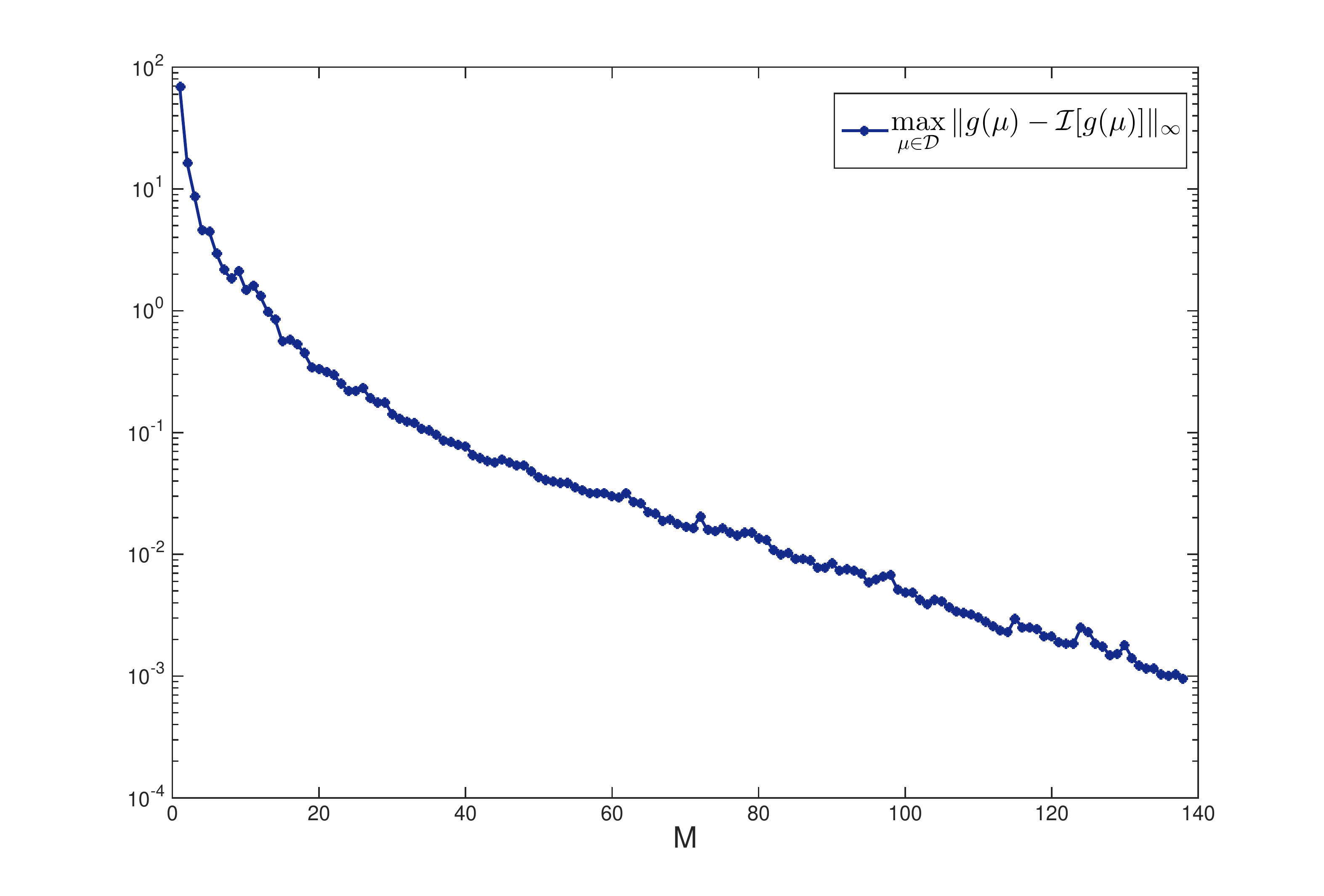}
\includegraphics[width=0.49\linewidth]{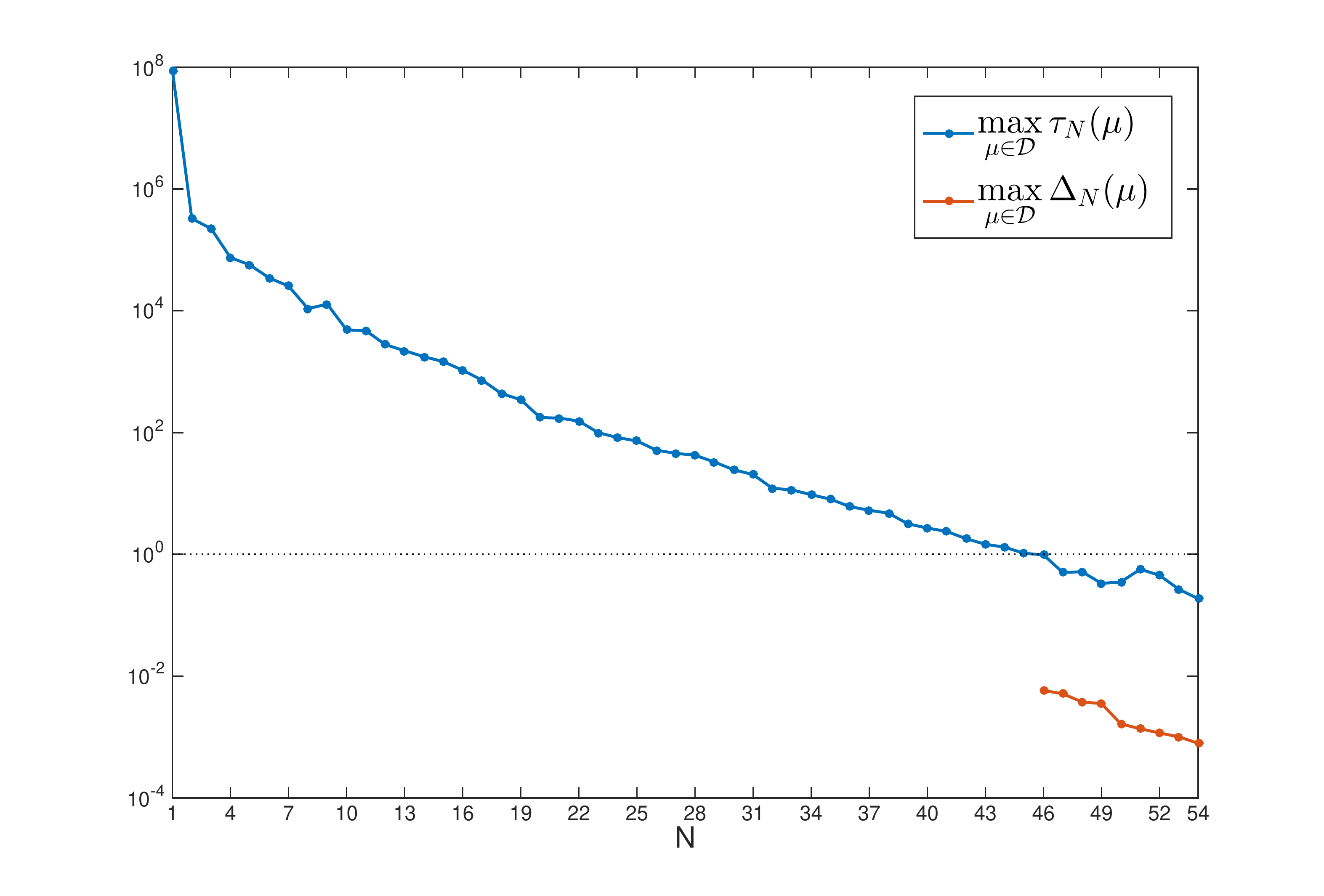}
\caption{Error evolution for the EIM, for Boussinesq VMS-Smagorinsky model with $\muk\in[10^3,10^4]\times[0.5,2]$.}\label{chap:Geom::fig:EIM_gRa}
\end{figure}

For the EIM, in this test, we prescribe a tolerance of $\veps_{EIM}=10^{-3}$. The error between $\nu_T(\uk_h;\muk)$ and its interpolant \black{fits} this tolerance when $M=138$ basis functions are included in the EIM reduced-basis space. In Fig. \ref{chap:Geom::fig:EIM_gRa} (left) we show the evolution of that error \black{along} the Greedy algorithm in the EIM.
  
\begin{figure}[h]
\centering
\includegraphics[width=0.75\linewidth]{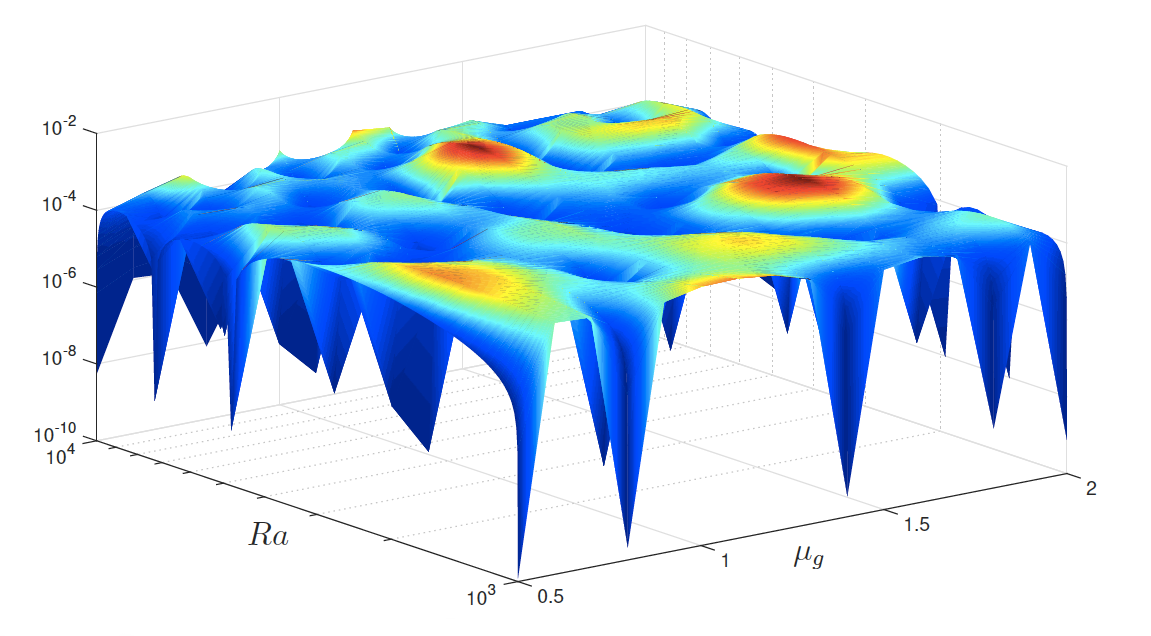}
\vspace{0.2cm}
\caption{\textit{A posteriori} error bound for $N_{\max}=54$.}\label{chap:Geom::fig:DeltaN_gRa}
\end{figure}

For the Greedy algorithm in the offline phase we prescribe a tolerance of $\veps_{RB}=10^{-3}$. This tolerance is reached when $N=N_{\max}=54$ basis functions are considered. We need $N=46$ basis functions to get $\tau_N(\muk)<1$, satisfying the conditions of Theorem \ref{chap:Geom::teor:Teorprinc}, and having defined the \textit{a posteriori} error bound $\Delta_N(\muk)$. In Fig. \ref{chap:Geom::fig:EIM_gRa} (left) we show the evolution of the maximum value of $\tau_N(\muk)$ and $\Delta_N(\muk)$ in the Greedy algorithm. On the other hand, in Fig. \ref{chap:Geom::fig:DeltaN_gRa} we show the value of the \textit{a posteriori} error estimator, when $N=N_{\max}=54$, for all $\muk\in\cD$.

%\begin{figure}[H]
%\centering
%\includegraphics[scale=0.3]{Est_err.eps}
%\caption{Evolution of the \textit{a posteriori} error bound in the Greedy algorithm, for Boussinesq VMS-Smagorinsky model with $\muk\in[10^3,10^4]\times[0.5,2]$.}\label{chap:Geom::fig:Delta_gRa}
%\end{figure} 

Finally, in Table \ref{chap:Geom::tab:geomRa} we sumarize some results obtained for some values of $\muk\in\cD$. There we show that the \red{relative error between the FE solution and the RB solution is between order $10^{-7}$ and $10^{-9}$ for velocity and temperature, and between order $10^{-8}$ and $10^{-9}$ for pressure}. For this test, the speedup rate obtained in the computation of the RB solution in the online phase with respect the computation of the FE solution is around fifty, due to the large number of EIM and RB basis functions needed, \red{ and the lower computational effort} in the FE computation. Again, we have obtained a good accuracy in the RB solution with respect to the FE solution, \black{with a considerable decrease} of the computational time. \black{For this test, the offline phase took approximately 2 weeks in being performed. In this offline computational time we consider either the EIM and the Greedy algorithm with the computation of the a posteriori error estimator.}

\begin{table}[H]
$$
\begin{tabular}{l|cccc}
\hline
Data &$Ra=2143$&$Ra=3506$&$Ra=5922$&$Ra=9618$\\
 &$\mu_g=1.95$&$\mu_g=0.71$&$\mug=1.13$&$\mug=1.63$\\
\hline
$T_{FE}$&600.96s & 914.18s & 684.95s&630.94s\\
$T_{online}$& 11.08s& 15.73s& 14.52s&11.46s\\
\hline
speedup& 54 & 58& 47 &55\\
\hline
$\|\uk_h-\uk_N\|_1/\|\uk_h\|_1$&$2.07\cdot10^{-7}$&$3.91\cdot10^{-7}$ & $6.52\cdot10^{-8}$ &$4.18\cdot10^{-8}$\\
\hline
$\|\theta_h-\theta_N\|_1/\|\theta_h\|_1$&$1.29\cdot10^{-7}$&$3.75\cdot10^{-7}$ & $3.17\cdot10^{-8}$ &$1.19\cdot10^{-8}$\\
\hline
$\|p_h-p_N\|_0\|p_h\|_0$&$1.56\cdot10^{-8}$&$4.21\cdot10^{-8}$ & $5.57\cdot10^{-9}$  &$1.52\cdot10^{-9}$\\
\hline
\end{tabular}$$
\hspace{-2cm}\caption{Computational time for FE and RB solutions, with the speedup and the error, for Boussinesq VMS-Smagorinsky model with $\muk\in[10^3,10^4]\times[0.5,2]$.}
\label{chap:Geom::tab:geomRa}
\end{table}

\section{Conclusions}\label{sec:Conclusions}

In this work, we have developed a reduced turbulence model for buoyant flows in domains with geometrical variability. \red{Specifically}, \black{we have dealt with} the RB Boussinesq VMS-Smagorinsky model, for a variable height cavity. To represent this variability in the cavity height, we have parametrized the domain. Thus, we needed to reformulate our problem in a reference domain, which \black{does not depend} on the geometric parameter. \black{As main technical tool,} we have developed an \textit{a posteriori} error estimator for the greedy algorithm involved in the \black{reduced basis space} construction. \black{This construction is based upon the Brezzi-Rappaz-Raviart theory. \black{We had to regularize} the eddy viscosity for temperature, in order to ensure that the Boussinesq-Smagorinsky operator is locally Lipschitz-continuous.}

Moreover, we have presented three different tests, considering \black{geometrical parameters, physical parameters, or both}. For each test, \black{we obtained an accurate} RB solution with a speedup rate going from  one thousand in the \black{simplest case}, to fifty in the most complex case \black{from one thousand for variability of only the physical parameter with diffusion-dominant effects, to nearly fifty for both geometrical and physical parameter variability}.

\section*{Acknowledgements}
This work has been supported by Spanish Government Project MTM2015-64577-C2-1-R and RTI2018-093521-B-C31, and COST Action TD1307 \\EU-MORNET. Francesco Ballarin and Gianluigi Rozza acknowledge project H2020 ERC CoG AROMA-CFD (GA 681447).

\section*{References}
\bibliography{Ref_RBM}

\appendix
\section{Proofs of theoretical results}\label{app:A}

\subsection{Proof of Proposition \ref{chap:Geom::prop::cont}}\label{app:prop1}
We consider $Z_h=V_h$ in $\partial_1A(U_h,V_h;\muk)(V_h)$. We first start bounding the diffusive terms for velocity and temperature, obtaining:
\begin{equation}\begin{array}{l}
a_{u,x}(\vk,\vk;\muk)+a_{u,y}(\vk,\vk;\muk)+a_{\theta,x}(\theta^v,\theta^v;\muk)
+a_{\theta,y}(\theta^v,\theta^v;\muk)\\
\ge \min\left\{\mug,\dfrac{1}{\mug}\right\}(Pr\normld{\nabla\vk_h}^2+\normld{\nabla\theta^v_h}^2).
\end{array}
\end{equation}
Denoting by $C_P$ the Poincare's constant, and considering the Holder's and Young's inequalities, the buoyancy term is bounded as
\begin{equation}
f(\theta^v,\vk;\muk)\ge -Pr\mug\mu\normld{\theta_h^v}\normld{\vk_h}
\ge -\dfrac{C_P Pr\mug\mu}{2}(\normld{\theta_h^v}^2+\normld{\vk_h}^2).
\end{equation}

Recalling the Sobolev embedding constants $C_u$ and $C_\theta$, defined in (\ref{chap:Geom::eq:Sobolev_u}) and (\ref{chap:Geom::eq:Sobolev_t}) respectively, we bound the following velocity convective terms as
\black{\begin{equation}
\begin{array}{c}
c_{u,x}(\vk_h,\uk_h,\vk_h;\muk)+c_{u,y}(\vk_h,\uk_h,\vk_h;\muk)\ge-|\min\{\mug,1\}(\vk_h\cdot\nabla\uk_h,\vk_h)_\Omega|\\
\ge-\min\{\mug,1\}\|\normlc{\vk_h}^2\normld{\uk_h}\ge -\min\{\mug,1\}C_u^2\normld{\nabla\uk_h}\normld{\nabla\vk_h}^2.
\end{array}
\end{equation}}

\black{The remaining convective terms can be bounded analogously, obtaining that }
\begin{equation}
\begin{array}{c}
c_{u,x}(\vk_h,\uk_h,\vk_h;\muk)+c_{u,x}(\uk_h,\vk_h,\vk_h;\muk)
+c_{u,y}(\vk_h,\uk_h,\vk_h;\muk)\\
+c_{u,y}(\uk_h,\vk_h,\vk_h;\muk)
+c_{\theta,x}(\vk_h,\theta_h^u,\theta_h^v;\muk)+c_{\theta,x}(\uk_h,\theta_h^v,\theta_h^v;\muk)\\
+c_{\theta,y}(\vk_h,\theta^u_h,\theta_h^v;\muk)+c_{\theta,y}(\uk_h,\theta_h^v,\theta_h^v;\muk)\\
\ge -\min\{\mug,1\}(2C_u^2\normld{\nabla\uk_h}\normld{\nabla\vk_h}^2
+2C_uC_\theta\normld{\nabla\theta_h^u}\normld{\nabla\theta^v}^2).
\end{array}
\end{equation}

For what concerns to the VMS-Smagorinsky terms, it holds
\begin{equation}
\begin{array}{c}
a_{Su,x}'(\uk_h;\vk_h,\vk_h;\muk) +a_{Su,y}'(\uk_h;\vk_h,\vk_h;\muk)
+a_{S\theta,nx}'(\uk_h;\theta^v_h,\theta^v_h;\muk)\\
+a_{S\theta,ny}'(\uk_h;\theta^v_h,\theta^v_h;\muk)\ge0,
\end{array}
\end{equation}
and
\begin{equation}
\begin{array}{c}
\mug\intOr{\partial_1\nu_T(\pih \uk)(\pih \zk)[\dx(\pih u_1)\dx(\pih v1)+\dx(\pih u_2)\dx(\pih v_2)]}\\
+\dfrac{1}{\mug}\intOr{\partial_1\nu_T(\pih \uk)(\pih \zk)[\dy(\pih u_1)\dy(\pih v1)+\dy(\pih u_2)\dy(\pih v_2)]} \\
\ge \min\left\{\mug,\dfrac{1}{\mug},1,\dfrac{1}{\mug^2}\right\}\intK{C_S^2\dfrac{1+\mug^2}{N_h^2}
\dfrac{|\nabla(\pih \uk_h):\nabla (\pih \vk_h)|^2}{|\nabla(T^{-1}\uk_h)|}}\ge0.
\end{array}
\end{equation}

Finally, using the local inverse inequalities (\textit{cf.} \cite{desinv}), we have that
\[
\dfrac{\mug}{Pr}\intOr{\partial_1\nu_{T,n}(\pih \uk)(\pih \zk)\;\dx(\pih \theta^u)\dx(\pih \theta^v)}
\]
\[
+\dfrac{1}{\mug\,Pr}\intOr{\partial_1\nu_{T,n}(\pih \uk)(\pih \zk)\;\dy(\pih \theta^u)\dy(\pih \theta^v})
\]
\[
\ge -C(\muk,\normlur{\phi_n'})
\normld{\nabla\uk_h}\normld{\nabla\theta^u_h}\normld{\nabla\vk_h}\normld{\nabla\theta^v_h}
\]
\[
\ge -\dfrac{C(\muk,\normlur{\phi_n'})}{2}\normld{\nabla\uk_h}\normld{\nabla\theta^u_h}(\normld{\nabla\vk_h}^2+\normld{\nabla\theta^v_h}^2),
\]
with
\[
C(\muk,\normlur{\phi_n'})=\dfrac{\min\left\{ \mug,\dfrac{1}{\mug} \right\} C_S^2\left(\sqrt{\dfrac{\mug^2+1}{N_h^2}}\right)^{2-d}C_f^4C\normlur{\phi_n'}}{Pr}.
\]

Thus, \black{taking into account all the previous bounds, we have proved that if (\ref{infsup::vel}) and (\ref{infsup::temp}) are verified, then there exists $\tilde{\beta}(\muk)>0$ such that}
\[
\partial_1A(U_h,V_h;\muk)(V_h)\ge \tilde{\beta}(\muk)(\normld{\nabla\vk_h}^2+\normld{\nabla\theta^v_h}^2) \forall V_h\in X_h.
\]
%\qed

\subsection{Proof of Proposition \ref{chap:Geom::prop:infsup}}\label{app:prop2}
We consider $Z_h=V_h$ in $\partial_1A(U_h,V_h;\muk)(V_h)$. We first start bounding the diffusive terms for velocity and temperature, obtaining:
\begin{equation}\begin{array}{l}
a_{u,x}(\vk,\vk;\muk)+a_{u,y}(\vk,\vk;\muk)+a_{\theta,x}(\theta^v,\theta^v;\muk)
+a_{\theta,y}(\theta^v,\theta^v;\muk)\\
\ge \min\left\{\mug,\dfrac{1}{\mug}\right\}(Pr\normld{\nabla\vk_h}^2+\normld{\nabla\theta^v_h}^2).
\end{array}
\end{equation}
Denoting by $C_P$ the Poincare's constant, and considering the Holder's and Young's inequalities, the buoyancy term is bounded as
\begin{equation}
f(\theta^v,\vk;\muk)\ge -Pr\mug\mu\normld{\theta_h^v}\normld{\vk_h}
\ge -\dfrac{C_P Pr\mug\mu}{2}(\normld{\theta_h^v}^2+\normld{\vk_h}^2).
\end{equation}

Recalling the Sobolev embedding constants $C_u$ and $C_\theta$, defined in (\ref{chap:Geom::eq:Sobolev_u}) and (\ref{chap:Geom::eq:Sobolev_t}) respectively, we bound the following velocity convective terms as
\black{\begin{equation}
\begin{array}{c}
c_{u,x}(\vk_h,\uk_h,\vk_h;\muk)+c_{u,y}(\vk_h,\uk_h,\vk_h;\muk)\ge-|\min\{\mug,1\}(\vk_h\cdot\nabla\uk_h,\vk_h)_\Omega|\\
\ge-\min\{\mug,1\}\|\normlc{\vk_h}^2\normld{\uk_h}\ge -\min\{\mug,1\}C_u^2\normld{\nabla\uk_h}\normld{\nabla\vk_h}^2.
\end{array}
\end{equation}}

\black{The remaining convective terms can be bounded analogously, obtaining that }
\begin{equation}
\begin{array}{c}
c_{u,x}(\vk_h,\uk_h,\vk_h;\muk)+c_{u,x}(\uk_h,\vk_h,\vk_h;\muk)
+c_{u,y}(\vk_h,\uk_h,\vk_h;\muk)\\
+c_{u,y}(\uk_h,\vk_h,\vk_h;\muk)
+c_{\theta,x}(\vk_h,\theta_h^u,\theta_h^v;\muk)+c_{\theta,x}(\uk_h,\theta_h^v,\theta_h^v;\muk)\\
+c_{\theta,y}(\vk_h,\theta^u_h,\theta_h^v;\muk)+c_{\theta,y}(\uk_h,\theta_h^v,\theta_h^v;\muk)\\
\ge -\min\{\mug,1\}(2C_u^2\normld{\nabla\uk_h}\normld{\nabla\vk_h}^2
+2C_uC_\theta\normld{\nabla\theta_h^u}\normld{\nabla\theta^v}^2).
\end{array}
\end{equation}

For what concerns to the VMS-Smagorinsky terms, it holds
\begin{equation}
\begin{array}{c}
a_{Su,x}'(\uk_h;\vk_h,\vk_h;\muk) +a_{Su,y}'(\uk_h;\vk_h,\vk_h;\muk)
+a_{S\theta,nx}'(\uk_h;\theta^v_h,\theta^v_h;\muk)\\
+a_{S\theta,ny}'(\uk_h;\theta^v_h,\theta^v_h;\muk)\ge0,
\end{array}
\end{equation}
and
\begin{equation}
\begin{array}{c}
\mug\intOr{\partial_1\nu_T(\pih \uk)(\pih \zk)[\dx(\pih u_1)\dx(\pih v1)+\dx(\pih u_2)\dx(\pih v_2)]}\\
+\dfrac{1}{\mug}\intOr{\partial_1\nu_T(\pih \uk)(\pih \zk)[\dy(\pih u_1)\dy(\pih v1)+\dy(\pih u_2)\dy(\pih v_2)]} \\
\ge \min\left\{\mug,\dfrac{1}{\mug},1,\dfrac{1}{\mug^2}\right\}\intK{C_S^2\dfrac{1+\mug^2}{N_h^2}
\dfrac{|\nabla(\pih \uk_h):\nabla (\pih \vk_h)|^2}{|\nabla(T^{-1}\uk_h)|}}\ge0.
\end{array}
\end{equation}

Finally, using the local inverse inequalities (\textit{cf.} \cite{desinv}), we have that
\[
\dfrac{\mug}{Pr}\intOr{\partial_1\nu_{T,n}(\pih \uk)(\pih \zk)\;\dx(\pih \theta^u)\dx(\pih \theta^v)}
\]
\[
+\dfrac{1}{\mug\,Pr}\intOr{\partial_1\nu_{T,n}(\pih \uk)(\pih \zk)\;\dy(\pih \theta^u)\dy(\pih \theta^v})
\]
\[
\ge -C(\muk,\normlur{\phi_n'})
\normld{\nabla\uk_h}\normld{\nabla\theta^u_h}\normld{\nabla\vk_h}\normld{\nabla\theta^v_h}
\]
\[
\ge -\dfrac{C(\muk,\normlur{\phi_n'})}{2}\normld{\nabla\uk_h}\normld{\nabla\theta^u_h}(\normld{\nabla\vk_h}^2+\normld{\nabla\theta^v_h}^2),
\]
with
\[
C(\muk,\normlur{\phi_n'})=\dfrac{\min\left\{ \mug,\dfrac{1}{\mug} \right\} C_S^2\left(\sqrt{\dfrac{\mug^2+1}{N_h^2}}\right)^{2-d}C_f^4C\normlur{\phi_n'}}{Pr}.
\]

Thus, \black{taking into account all the previous bounds, we have proved that if (\ref{infsup::vel}) and (\ref{infsup::temp}) are verified, then there exists $\tilde{\beta}(\muk)>0$ such that}
\[
\partial_1A(U_h,V_h;\muk)(V_h)\ge \tilde{\beta}(\muk)(\normld{\nabla\vk_h}^2+\normld{\nabla\theta^v_h}^2) \forall V_h\in X_h.
\]
%\qed

\subsection{Proof of Lemma \ref{chap:Geom::lema:LemmaRho}}\label{app:Lemma}
Thanks to the triangular inequality, it holds
\[
\left|\partial_1A(U_h^1,V_h;\muk)(Z_h)-\partial_1A(U_h^2,V_h;\muk)(Z_h)\right|
\]
\[
\le\max\{1,\mug\}|(\zk_h\cdot\nabla(\uk^1_h-\uk^2_h),\vk_h)_{\Or}|+\max\{1,\mug\}|((\uk^1_h-\uk^2_h)\cdot\nabla\zk_h,\vk_h)_{\Or}|
\]
\[
+\max\{1,\mug\}
|(\zk_h\cdot\nabla(\theta^{u1}_h-\theta^{u2}_h),\theta^v_h)_{\Or}|+\max\{1,\mug\}|((\uk^1_h-\uk^2_h)\cdot\nabla\theta^z_h,\theta^v_h)_{\Or}|
\]
\[
+\max\left\{\mug,\dfrac{1}{\mug}\right\}|(\nu_T(\uk_h^1;\muk)-\nu_T(\uk_h^2;\muk))\nabla(\pih\zk_h),\nabla(\pih\vk_h))_{\Or}|
\]
\[
+\max\left\{\mug,\dfrac{1}{\mug}\right\}|(\nu_{T,n}(\uk_h^1;\muk)-\nu_{T,n}(\uk_h^2;\muk))\nabla(\pih\theta^z_h),\nabla(\pih\theta^v_h))_{\Or}|
\]
\[
+\max\left\{\mug,\dfrac{1}{\mug}\right\}\Big|(\partial_1\nu_T(\pih\uk^1_h)(\pih\zk_h)\nabla(\pih\uk_h^1),\nabla(\pih\vk_h))_{\Or}
\]
\[
-(\partial_1\nu_T(\pih\uk^2_h)(\pih\zk_h)\nabla(\pih\uk_h^2),\nabla(\pih\vk_h))_{\Or}\Big|
\]
\[
+\max\left\{\mug,\dfrac{1}{\mug}\right\}\Big|(\partial_1\nu_{T,n}(\pih\uk^1_h)(\pih\zk_h)\nabla(\pih\theta^{u1}_h),\nabla(\pih\theta^v_h))_{\Or}
\]
\[
-(\partial_1\nu_{T,n}(\pih\uk^2_h)(\pih\zk_h)\nabla(\pih\theta_h^{u2}),\nabla(\pih\theta^v_h))_{\Or}\Big|
\]

We bound each term separately. The first four terms, corresponding with the convective terms are bounded analogously. For brevity, we show one of them. Thus, considering the Sobolev embedding constants (\ref{chap:Geom::eq:Sobolev_u}) and (\ref{chap:Geom::eq:Sobolev_t}),
\begin{equation}
\begin{array}{c}
|(\zk_h\cdot\nabla(\theta^{u1}_h-\theta^{u2}_h),\theta^v_h)_{\Or}|\le \normlc{\zk_h}\normld{\nabla(\theta^{u1}_h-\theta^{u2}_h)}\normlc{\theta^v_h}\vspace{0.1cm}\\ 
\le C_uC_{\theta}\normX{U_h^1-U_h^2}\normX{Z_h}\normX{V_h}
\end{array}
\end{equation}

The VMS-Smagorinsky terms for eddy viscosity and eddy diffusivity are also bounded in a similar way. We show the boundness of the eddy diffusivity term, for which we take into account inequality (\ref{eq:pih_ineq}), the inverse inequalities (\textit{cf.} \cite{desinv}) and the properties of the convolution, recalling that $\normlur{\phi_n}=1$,
\begin{equation}
\begin{array}{c}
|(\nu_{T,n}(\uk_h^1;\muk)-\nu_{T,n}(\uk_h^2;\muk))\nabla(\pih\theta^z_h),\nabla(\pih\theta^v_h))_{\Or}|\vspace{0.1cm}\\
\le (C_Sh)^2\normlur{\phi_n}\normlt{\nabla(\pih \uk^1_h-\pih \uk_h^2)}\normlt{\nabla(\pih\theta^z_h)}\normlt{\nabla(\pih\theta^v_h)}\vspace{0.1cm}\\
\le C_S^2h^{2-d/2}C_f^3\normX{U_h^1-U_h^2}\normX{Z_h}\normX{V_h}
\end{array}
\end{equation}

Taking into account \black{inequality (\ref{eq:pih_ineq}), the seventh} term can be bounded as in \red{Lemma} 5.1 of \cite{PaperSmago}. \black{We resume the bound in the following
\[
\Big|(\partial_1\nu_T(\pih\uk^1_h)(\pih\zk_h)\nabla(\pih\uk_h^1),\nabla(\pih\vk_h))_{\Or}\]\[
-(\partial_1\nu_T(\pih\uk^2_h)(\pih\zk_h)\nabla(\pih\uk_h^2),\nabla(\pih\vk_h))_{\Or}\Big|
\]
\[
\le (C_Sh)^2\normlt{\nabla\zk}\normlt{\nabla(\uk_h^1-\uk_h^2)}\normlt{\nabla\vk_h}
\]
\[
+(C_Sh)^2\normlt{\nabla(\uk_h^1-\uk_h^2)}\normlt{\nabla\zk}\normlt{\nabla\vk_h}
\]
\[
+(C_Sh)^2\normlt{\nabla(\uk_h^1-\uk_h^2)}\normlt{\nabla\zk}\normlt{\nabla\vk_h}
\]
\[
\le 3C_Sh^{2-d/2}C\normld{\nabla(\uk_h^1-\uk_h^2)}\normld{\nabla\zk}\normld{\nabla\vk_h}
\]
\[
\le 3C_Sh^{2-d/2}C\normX{U_h^1-U_h^2}\normX{Z_h}\normX{V_h}.
\]
}

Finally, we next show the bound of the last term, taking into account again the inverse inequalities and the Sobolev embedding constants: 
\[
\big|(\partial_1\nu_{T,n}(\pih\uk^1_h)(\pih\zk_h)\nabla(\pih\theta^{u1}_h),\nabla(\pih\theta^v_h))_{\Or}\]\[
-(\partial_1\nu_{T,n}(\pih\uk^2_h)(\pih\zk_h)\nabla(\pih\theta^{u2}_h),\nabla(\pih\theta^v_h))_{\Or}\big|
\]\[
\le \big|(\partial_1\nu_{T,n}(\pih\uk^1_h)(\pih\zk_h)\nabla(\pih(\theta^{u1}_h-\theta^{u2}_h),\nabla(\pih\theta^v_h))_{\Or}\big|
\]\[
+\big|[(\partial_1\nu_{T,n}(\pih\uk^2_h)(\pih\zk_h)-(\partial_1\nu_{T,n}(\pih\uk^2_h)(\pih\zk_h)]\nabla(\pih\theta^{u2}_h),\nabla(\pih\theta^v_h))_{\Or}\big|
\]\[
\le C_1(\normlinf{\pih(\nabla\uk^1_h)})\normld{\nabla\zk_h}\normld{\nabla(\theta^{u1}_h-\theta^{u2}_h)}\normld{\theta_h^v}
\]\[
+C_2(\normlinf{\pih(\nabla\theta^{u2}_h)})\normld{\nabla\zk_h}\normld{\nabla(\uk_h^1-\uk_h2)}\normld{\theta_h^v}
\]\[
\le [C_1(\normlinf{\pih(\nabla\uk^1_h)}) +C_2(\normlinf{\pih(\nabla\theta^{u2}_h)})]\normX{U_h^1-U_h^2}\normX{Z_h}\normX{V_h}.
\]
%\qed

\subsection{Proof of Theorem \ref{chap:Geom::teor:Teorprinc}}\label{app:Teorema}
This proof is an adaptation of the proofs of Theorems 5.2 and 5.3 of \cite{PaperSmago} and Theorem 3.3 of \cite{Deparis}.

We  define the following operators:
\begin{itemize}
\item $\cR{\cdot}:X_h\ra X'_h$, defined as 
\begin{equation}\label{cR}
\left<\cR{Z_h},V_h\right>=A(Z_h,V_h;\muk)-F(V_h;\muk), \quad\forall Z_h,V_h\in X_h
\end{equation}
\item $\DA{U_h(\muk)} :X_h\ra X'_h$, defined, for $U_h(\mu)\in X_h$, as 
\begin{equation}\label{Thmu}
\left<\DA{U_h(\muk)}Z_h,V_h\right>=\partial_1A(U_h(\muk),V_h;\muk)(Z_h),\quad\forall Z_h,V_h \in X_h
\end{equation}
\item $H: X_h\ra X_h$, defined as
\begin{equation}\label{H}
H(Z_h;\muk)=Z_h-\DAn^{-1}\cR{Z_h}, \quad\forall Z_h\in X_h
\end{equation}
\end{itemize}

Note that $\DAn$ is invertible thanks to the assumption $\beta_N(\mu)>0$. . %To prove the existence and uniqueness of problem ($\ref{pb}$), we have to prove that $H$ is a contractive mapping on $X_h$ and then it admits a unique fixed point. For this purpose, 
We express
\begin{equation}\label{difH}
\difH=(\difZ)-\DAn^{-1}(\cR{Z_h^1}-\cR{Z_h^2}).
\end{equation}

It holds
\begin{equation}\label{difRes}
\cR{Z_h^1}-\cR{Z_h^2}=\DA{\xi}(\difZ),
\end{equation}
where $\xi=\lambda Z_h^1-(1-\lambda)Z_h^2,$ for some $\lambda\in(0,1)$.  Multiplying (\ref{difH}) by $\DAn$ and applying this last property, we can write
\[
\DAn(\difH)=\left[\DAn-\DA{\xi}\right](\difZ).
\]

Then, thanks to Lemma \ref{chap:Geom::lema:LemmaRho} and this last equality, it follows that in a neighborhood of $\unmu$ and $\xi$,
\[
\left<\DAn(\difH),V_h\right>
\]\[
\le\rho_n(\mug)\normX{\unmu-\xi}\normX{\difZ}\|V_h\|_X.
\]

Now, applying the definitions of $\beta_N(\mu)$, $T_N$, $\DAn$, and this last property, we can obtain
\[
\beta_N(\mu)\normX{\difH}\|T_N(\difH)\|_X
\]
\[
\le\|T_N(\difH)\|_X^2
\]\[
=\big(\TN{\difH},\TN{\difH}\big)_X
\]\[
=\left<\DAn(\difH), \TN{\difH}\right>
\]\[
\le\rho_n(\mug)\normX{\unmu-\xi}\normX{\difZ}\|\TN{\difH}\|_X%=\romu\normu{\unmu-\xi}\normu{\difu}.
\]

So, we have proved that
\[
\normX{\difH}\le\frac{\rho_n(\mug)}{\beta_N(\muk)}\normX{\unmu-\xi}\normX{\difZ}.
\]

If $Z_h^1$ and $Z_h^2$ are in $B_X(\unmu,\alpha)$ then, $\normX{\unmu-\xi}\le\alpha$, and,
\[
\normX{\difH}\le\frac{\rho_n(\mug)}{\beta_N(\muk)}\alpha\normX{\difZ}.
\]

Then, $\Hw{\cdot}$ is a contraction if $\alpha<\dfrac{\beta_N(\muk)}{\rho_n(\mug)}.$ So it follows that there can exist at most one fixed point of $\Hw{\cdot}$ inside $B_X\left(\unmu,\dfrac{\beta_N(\muk)}{\rho_n(\mug)}\right)$, and hence, at most one solution $U_h(\muk)$ to (\ref{chap:Geom::pb:FVX}) in this ball.

To prove (\ref{chap:Geom::teorprinc:err}), we prove that the operator $\Hw{\cdot}$ has a fixed point.
Thus, let $\alpha>0$ and $Z_h\in X_h$ such that $\normX{\unmu-Z_h}\le\alpha$. We consider
\[
\Hw{Z_h}-\unmu=Z_h-\unmu-\DAn^{-1}\cR{Z_h}
\]\[
=Z_h-\unmu-\DAn^{-1}\left[\cR{Z_h}-\cR{\unmu}\right]
\]\[
-\DAn^{-1}\cR{\unmu}
\]

Multiplying by $\DAn$, we obtain
\[
\left<\DAn(\Hw{Z_h}-\unmu),V_h\right>=\left<\DAn(Z_h-\unmu),V_h\right>
\]\[
-\left<\cR{Z_h}-\cR{\unmu},V_h\right>-\left<\cR{\unmu},V_h\right>,\quad\forall V_h\in X_h.
\]

It holds that $\cR{Z_h}-\cR{\unmu}=\DA{\xi(\mu)}(Z_h-\unmu)$, where $\xi(\mu)=t^*Z_h+(1-t^*)\unmu$, $t^*\in(0,1)$.

Thus, Lemma \ref{chap:Geom::lema:LemmaRho} and this last equality, it follows that in a neighborhood of $\unmu$ and $\ximu$, we obtain:
\[
\left<\DAn(\Hw{Z_h}-\unmu),V_h\right>=\dual{\DAn(Z_h-\unmu)}
\]\[
-\dual{\DA{\ximu}(Z_h-\unmu)}-\dual{\cR{\unmu}}
%%\]\[
%%=\dual{\big(\DAn-\DA{\ximu}\big)(Z_h-\unmu)}-\dual{\cR{\unmu}}
\]\[
\le\rho_n(\mug)\normX{\unmu-\ximu}\normX{Z_h-\unmu}\|V_h\|_X+\en\|V_h\|_X
\]\[
\le\big(\rho_n(\mug)\normX{Z_h-\unmu}^2+\en\big)\normX{V_h}
\]

Thus, it follows that,
\[
\beta_N(\mu)\normX{\Hw{Z_h}-\unmu}\normX{\Tnmu(\Hw{Z_h}-\unmu)}
\]
\[
\le\|\Tnmu(\Hw{Z_h}-\unmu)\|_{X}^2
%\]\[
%=\Big(\Tnmu\big(\Hw{Z_h}-\unmu\big),\,\Tnmu\big(\Hw{Z_h}-\unmu\big)\Big)_X
%\]\[
%=\left< \DAn(\Hw{Z_h}-\unmu),\Tnmu(\Hw{Z_h}-\unmu)\right>
\]\[
\le\big(\rho_n(\mug)\normX{Z_h-\unmu}^2+\en\big)\normX{\Tnmu(\Hw{Z_h}-\unmu)}.
\]

Then, as $Z_h\in B_X\left(\unmu,\alpha\right)$, we have 
\begin{equation}\label{seggrado}
\normX{\Hw{Z_h}-\unmu}<\dfrac{\rho_n(\mug)}{\bmu}\alpha^2+\dfrac{\en}{\bmu}.
\end{equation}

In order to ensure that $H$ maps $B_X(\unmu,\alpha)$ into a part of itself, we are seeking the values of $\alpha$ such that $
\dfrac{\rho_n(\mug)}{\bmu}\alpha^2+\dfrac{\en}{\bmu}\le\alpha$. This condition is verified for $\alpha=\dnk$. Consequently, since $\dnk\le \dfrac{\beta_N(\muk)}{\rho_n(\mug)}$, there exists a unique solution $U_h(\mu)$ to (\ref{chap:Geom::pb:FVX}) in the ball  $B_X(\unmu,\alpha)$.

Finally, (\ref{chap:Geom::teorprinc:efec}), can be proved analogously as in \cite{PaperSmago}.\\
%\qed

\end{document}